\documentstyle[12pt]{amsart}
\newtheorem{th}{Theorem}[section]
\newtheorem{crl}[th]{Corollary}
\newtheorem{prp}[th]{Proposition}
\newtheorem{lm}[th]{Lemma}

\newcommand{\der }{\partial }

\newcommand{\eps}{\epsilon}

\newcommand{\limline}{\limits_{-\!-\!-\!-\!-\!-}}
\newcommand{\limsim}{\limits_{\sim\!\sim\!\sim\!\sim\!\sim\!\sim}}
\newcommand{\limeq}{\limits_{=\!=\!=\!=\!=\!=}}
\newcommand{\limsimeq}{\limits_{\simeq\!\simeq\!\simeq\!\simeq\!\simeq\!\simeq}}
\newcommand{\limequiv}{\limits_{\equiv\!\equiv\!\equiv\!\equiv\!\equiv\!\equiv}}
\newcommand{\limsmile}{\limits_{\smile\!\smile\!\smile\!\smile\!\smile}}
\newcommand{\limapprox}{\limits_{\approx\!\approx\!\approx\!\approx\!\approx\!\approx}}
\newcommand{\limcong}{\limits_{\cong\!\cong\!\cong\!\cong\!\cong\!\cong}}
\newcommand{\limasymp}{\limits_{\asymp\!\asymp\!\asymp\!\asymp\!\asymp\!\asymp}}
\newcommand{\limdoteq}{\limits_{\doteq\!\doteq\!\doteq\!\doteq\!\doteq\!\doteq}}
\newcommand{\limfrown}{\limits_{\frown\!\frown\!\frown\!\frown\!\frown}}
\newcommand{\limequi}{\limits_{\equiv\; \equiv\; \equiv\; \equiv}}
\newcommand{\limcon}{\limits_{\cong\; \cong\; \cong\; \cong}}
\newcommand{\limsi}{\limits_{\sim\; \sim\; \sim\; \sim}}
\newcommand{\limli}{\limits_{-\; -\; -\; -}}
\newcommand{\lime}{\limits_{=\; =\; =\; =}}
\newcommand{\limsime}{\limits_{\simeq\; \simeq\; \simeq\; \simeq}}

\begin{document}
\title[Cohomologies and deformations of 
right-symmetric \ldots]{Cohomologies and
deformations of right-symmetric algebras}
\author {Askar Dzhumadil'daev}
\address
{Inst. Mathematics, Pushkin str.125, Almaty, 480021,
Kazakhstan}
\email{askar@@itpm.sci.kz}
\maketitle

{\em Dedicated to A.I. Kostrikin on the occasion of his 70-th birthday}

\begin{abstract}
An algebra  $A$ with identity 
$(a\circ b)\circ c-a\circ(b\circ c)=(a\circ c)\circ b-a\circ(c\circ b),$
is called right-symmetric. Cohomology and deformation theory for
right-symmetric algebras are developed. 
Cohomologies of $gl_n$ and half-Witt algebras $W_n^{rsym},\;  p=0,$ 
$W_n^{rsym}({\bf m}),\; p>0,$ are
calculated. In particular, one right-symmetric central extension of 
$W_1^{rsym}$ is constructed.

\end{abstract}

\tableofcontents

\section{Introduction.}
An algebra $A$ over a field ${\cal K}$ of characteristic $p\ge 0$ 
is called right-symmetric \cite{Vinberg}, \cite{Koszul}, if for any
$a,b,c\in A,$ the following condition takes place
$$a\circ(b\circ c)-(a\circ b)\circ c= a\circ (c\circ b)-(a\circ c)\circ b.$$
Any associative algebra is right-symmetric. For example, $gl_n$ under
usual multiplication of matrices is right-symmetric.
An algebra of vector fields ${\cal K}[[x^{\pm 1},\ldots, x^{\pm 1}_n]]$ under multiplication
$u\der_i\circ v\der_j=v\der_j(u)\der_i$ gives us a less trivial 
example of right-symmetric algebras. It is not associative.
 Since its Lie algebra is isomorphic to
Witt algebra $W_n,$ we call it a half-Witt algebra and denote it as 
$W_n^{rsym}.$ If $n=1,$ this algebra satisfies one more identity
$$a\circ (b\circ c)=b\circ (a\circ c).$$
Such algebras are called Novikov \cite{Nov}, \cite{Osborn1}. 
The generalisation of Novikov structure for the case $n>1$ is possible, 
if we consider half-Witt algebra not with one, 
but with two multiplications. If we endow ${\cal K}[[x^{\pm 1},\ldots,
x^{\pm 1}]]$ with a second multiplication $u\der_i\ast v\der_j= 
\der_i(u)v\der_j,$ then we obtain algebra with the following identities
$$a\circ (b\circ c)-(a\circ b)\circ c-a\circ(c\circ b)+(a\circ c)\circ b=0,$$
$$a\ast (b\ast c)-b\ast(a\ast c)=0,$$
$$a\circ (b\ast c)-b\ast(a\circ c)=0,$$
$$(a\ast b-b\ast a- a\circ b+b\circ a)\ast c=0,$$
$$(a\circ b- b\circ a)\ast c+a\ast(c\circ b)-(a\ast c)\circ b
-b\ast(c\circ a)+(b\ast c)\circ a=0.$$

\noindent These two multiplications are useful in the construction of 
right-symmetric, Chevalley-Eilenberg and Leibniz cocycles of such algebras.

We develop cohomology theory for right-symmetric algebras. 
We endow right-symmetric cochain complex $C^*_{rsym}(A,M)=
\oplus_kC^k_{rsym}(A,M),$ where $C^{k+1}_{rsym}(A,M)=Hom(A\otimes
\wedge^k(A),M), k\ge 0,$ by a pre-simplicial structure. 
Corresponding cohomologies can be "almost" obtained by derived functor 
formalism. The exact meaning of the word  "almost" can be found in section
\ref{derived}. Roughly speaking, this means that one should be more careful 
in considerating small degree cohomologies. If we take $C^0_{rsym}(A,M)$ 
as $M,$ then we should consider the operator $d_{rsym}$ with cubic condition 
$d_{rsym}^3=0.$ \cite{Kerner}. We prefer taking $C^0_{rsym}(A,M)$ as a 
$Ker\,d^2_{rsym}$ on $M,$ i.e., $C^0_{rsym}(A,M):=M^{l.ass}:=\{m\in M: 
(m,a,b)=0, \forall a,b\in A\}.$ Then for any $m\in M,$ we can correspond 
2-right-symmetric cocycles, $\nabla(m):(a,b)\mapsto (m,a,b).$ We call such 
cocycles as standard. If $m\in M^{l.ass},$ then $\nabla(m)=0.$ 
If $m\in M,$ then the cohomological class $[\nabla(m)]=0,$ because of
$\nabla(m)=d\omega,$ where $\omega(a)=[a,m].$ Moreover it is true, 
if $M$ is a submodule of some right-symmetric $A-$module $\tilde M,$ 
and $m\in \tilde M,$ such that $[a,m]=a\circ m-m\circ a\in M, 
\forall a\in M.$ If $\tilde m\in \tilde M,$ such that $d_{rsym}\tilde 
m(a)\not\in M,$ then $\nabla({\tilde m})$ can give a nontrivial 
class of 2-right-symmetric cocycles in $H^2_{rsym}(A,M).$ For example,
Osborn 2-right-symmetric cocycles for $A=W_1^{rsym}(m), p>0,$ that appear
in constructing simple Novikov algebras,
$$(u\der,v\der)\mapsto x^{p^m-1}uv\der,$$
$$(u\der,v\der)\mapsto x^{p^m-2}uv\der,$$

\noindent are $\nabla({x^{p^m+1}\der}),$ and $\nabla({x^{p^m}\der}),$ 
correspondingly.

If $k>0,$ right-symmetric cohomologies $H^{k+1}_{rsym}(A,M)$
are isomorphic to Chevalley-Eilenberg cohomologies
$H^k_{lie}(A,C^1(A,M)),$ where $A^{lie}-$module structure on $C^1(A,M)$
is given by a special way: $[a,f](b)= -d_{rsym}f(b,a).$ 
We endow also right-symmetric universal enveloping algebra by a 
Hopf algebraic structure. It allows us to consider cup products,
that are very useful in cocycle constructions.

Second cohomology space $H^2_{rsym}(A,A)$ is interpreted as a space 
of right-symmetric deformations. We calculate right-symmetric cohomologies of 
matrix algebra $gl_n, p=0.$ We prove that, in the category of irreducible
antisymmetric $gl_n^{rsym}-$modules, nontrivial cohomologies appear only in 
the case of $M=(gl_n)_{anti}.$  Moreover, right-symmetric cohomology of 
$gl_n^{rsym}$ in $(gl_n)_{anti}$ can be reduced to Chevalley-Eilenberg cohomology 
of Lie algebra $gl_n$ with coefficients in trivial module:
$$H^{k+1}_{rsym}(gl_n,(gl_n)_{anti})\cong H^k_{lie}(gl_n,{\cal K}),\; k>0.$$
In particular, $H^{k+1}_{rsym}(gl_n,{\cal K})=0, \; k\ge 0.$
We calculate also right-symmetric cohomologies of $gl_n$ with coefficients in 
regular module. These results show that $gl_n$ has $(n^2-1)-$parametrical 
nontrivial right-symmetric deformations. Any formal right-symmetric deformation of
$gl_n$ is equivalent to the deformations given by the rule
$$(a,b)\mapsto a\circ b + t\,tr\,b\,[X,a],\;\; X\in sl_2.$$

\noindent One can choose the  prolongation in another way: 
$$(a,b)\mapsto a\circ b+ t\,X\circ ((tr\,a)b-tr\,(a\circ b) +(tr\,b)a)$$
$$+t^2\{tr\,a\,tr\,b-(tr \,a\circ b)^2X^2-(tr\,a\,tr(X\circ b))X-(tr\,(a\circ X)
tr\,b)X\}+\cdots .$$

We prove that right-symmetric cohomologies of $A=W_n^{rsym}$  with
coefficients in antisymmetric modules can also be reduced to 
Chevalley-Eilenberg cohomologies of the Lie algebra $W_n.$
As it turned out, $H^2_{rsym}(A,A)$ for $A=W_n, p=0,$ or $A= W_n({\bf m}),
p>0,$ is too large and this happens mainly because of largeness of a space
of right-symmetric derivations. There is an imbedding 
$$Z^1_{rsym}(A,A)\otimes H^1_{lie}(A,U)\rightarrow H^2_{rsym}(A,A).$$

\noindent We prove that $Z^1_{rsym}(A,A)$ has a basis consisting of two 
types of right-symmetric derivations: $\der_i, i=1,\ldots,n,$ if $p>0,$ 
one should consider also derivations 
$\der_i^{p^{k_i}}, 0\le k_i<m_i;$ and $x_i\der_j, i,j=1,\ldots,n.$ So,
any right-symmetric derivation of $A$ has a form 
$\sum_{i=1}^nu_i\der_i+\delta(p>0)\sum_{i=1}^n\sum_{k_i=0}^{m_i-1}
\lambda_{i,k_i}\der_i^{p^{k_i}},$ such that $\der_i\der_j(u_s)=0,
i,j,s=1,\ldots,n, \lambda_{i,k_i}\in {\cal K}.$
We formulate a result about local deformations of $W_n,  p=0,$ or
$W_n({\bf m}), p>3.$ 
The space $H^2_{rsym}(W_n,W_n), p=0,$ is generated by classes of cocycles
of  four types. In the case of $p>3$ Steenrod Squares also appears. 
We prove that $W_1^{rsym}$ has exactly one right-symmetric central extension.
It can be given by cocycle
$$(e_i,e_j)\mapsto (j+1)j\delta_{i+j,-1},\;\; p=0,$$
$$(e_i,e_j)\mapsto (-1)^i\delta_{i+j,p^m-1},\;\; p>0.$$
For $n>1,$ $H^2_{rsym}(W_n,{\cal K})=0.$

For right-symmetric algebras, Novikov algebras and some cohomology
calculations see also \cite{GelfandDorfman}, \cite{Kim}, \cite{Osborn2}, 
\cite{Osborn3}, \cite{Burde1}, \cite{Burde2}, \cite{Shen}.

\section{Right-symmetric algebras and (co)modules.\label{1}}

\subsection{Right-symmetric algebras}
An algebra $A$ over a field ${\cal K}$ with multiplication
$(a,b)\mapsto a\circ b,$ is called {\it Lie-admissible,} if the vector
space $A$  under commutator $[a,b]= a\circ b-b\circ a$ can be endowed by
a structure of Lie algebra. An algebra $A$ is {\it right-symmetric,} if it
satisfies the following identity :
$$(a\circ b)\circ c-a\circ(b\circ c)=(a\circ c)\circ b-a\circ(c\circ b), \quad
\forall a,b,c\in A.$$

\noindent Let $(a,b,c)=a\circ(b\circ c)-(a\circ b)\circ c$ be the asssociator
of elements $a,b,c\in A.$ In terms of associators the right-symmetric identity
is
$$(a,b,c)=(a,c,b), \; \forall a,b,c\in A.$$

\noindent Right-symmetric algebra $A$ is Lie-admissible. Similarly, one can define
left-symmetric algebra by identity
$$(a,b,c)=(b,a,c), \; \forall a,b,c\in A.$$

\noindent Categories of left-symmetric algebras and right-symmetric algebras are
equivalent. Any left(right)-symmetric algebra will be right(left)-symmetric under
new multiplication $(a,b)\mapsto b\circ a.$

An element $e$ of right-symmetric algebra is called {\it left unit}, if
$e\circ a=a,$ for any $a\in A.$ Denote by $Q_l(A)$ a space of left units.
Let $Z_l(A)=\{z\in A: z\circ a=0,\forall a\in A\}$ be {\it left center}
of  $A.$ Call a space $N_l(A)= Z_l(A)\oplus Q_l(A)$ as a semi-senter of
$A.$ Then $[N_l(A),N_l(A)]\subseteq Z_l(A).$  An algebra $A$ is called
{\it (left)} unital, if it has nontrivial left units.

Any associative algebra is a right-symmetric algebra.
In such cases, we will use notations like $A^{ass},$ if we consider
$A$ as associative algebra and $A^{rsym},$ if we consider $A$ as
right-symmetric algebra. Similarly, for right-symmetric algebra
$A$ notation $A^{rsym}$ means that we use only right-symmetric 
structure on $A$ and $A^{lie}$ stands for a Lie algebra structure 
under commutator $(a,b)\mapsto [a,b].$

Matrix algebras $gl_n$ gives us examples of unital right-symmetric algebras.

Less trivial examples appear in the consideration of Witt algebras. The algebra
$W_n, p=0, $ and $W_n({\bf m})$ defined below has not only right-symmetric
mutiplication $(a,b)\mapsto a\circ b,$ but also one more
multiplication $(a,b)\mapsto a\ast b,$ that satsifies the following
identities
$$a\circ (b\circ c)-(a\circ b)\circ c-a\circ(c\circ b)+(a\circ c)\circ b=0,$$
$$a\ast (b\ast c)-b\ast(a\ast c)=0,$$
$$a\circ (b\ast c)-b\ast(a\circ c)=0,$$
$$(a\ast b-b\ast a- a\circ b+b\circ a)\ast c=0,$$
$$(a\circ b- b\circ a)\ast c+a\ast(c\circ b)-(a\ast c)\circ b
-b\ast(c\circ a)+(b\ast c)\circ a=0.$$

Let
$$U=k[[x^{\pm 1}_1,\ldots,x^{\pm 1}_n]]=\{x^{\alpha}=\prod_{i=1}^kx_i^{\alpha_i}:
\alpha=(\alpha_1,\ldots,\alpha_n), \alpha_i\in {\bf Z}, i=1,\ldots,n\}$$

\noindent be an algebra of Laurent power series, if the main field $k$ has
characteristic 0 and
$$U=O_n({\bf m})=\{x^{(\alpha)}=\prod_{i}x_i^{(\alpha_i)}:
\alpha=(\alpha_1,\ldots,\alpha_n), 0\le \alpha_i< p^{m_i}, i=1,\ldots,n\}$$

\noindent be a divided power algebra if $char\, k=p>0.$ Recall that  $O_n({\bf m})$ is
$p^m$~-dimensional and the muliplication is given by
$$x^{(\alpha)}x^{(\beta)}={{\alpha+\beta}\choose \alpha}
x^{(\alpha+\beta)},$$

\noindent where $m=\sum_im_i,$ and
$${{\alpha+\beta}\choose \alpha}=\prod_{i}{{\alpha_i+\beta_i}\choose \alpha_i},
\; {n \choose l}= {n!\over {l!(n-l)!}}, \; n,l\in {\bf Z}_+.$$

Let $\eps_i=(0,\ldots,\mathop{1}\limits_{i},\ldots,0).$ Define
$\der_i$ as a derivation of $U,$
$$\der_i(x^\alpha)=\alpha_ix^{\alpha-\eps_i}, \; p=0,$$
$$\der_i(x^{(\alpha)})=x^{(\alpha-\eps_i)}, \; p>0.$$

\noindent Endow a space of derivations $Der\,U=\{\sum_iu_i\der_i: u_i\in U\}$
by multiplications:
$$u\der_i\circ v\der_j=v\der_j(u)\der_i,$$
$$u\der_i\ast v\der_j=\der_i(u)v\der_j.$$

\noindent Denote obtained algebra as $W_n^{rsym}({\bf m}).$
If $p=0,$ this denotion will be reduced until $W_n^{rsym}.$
If $n\not\equiv 0(mod\,p),$ then an element $e=\sum_ix_1\der_i/n$ is
a left unit of $W_n^{rsym}(m).$ If $n=1,$ then $a\circ b=a\ast b.$  
Thus the algebra $A=W_1^{rsym}({\bf m})$ in addition to right-symmetry
condition satisfies the following identity
$$a\circ (b\circ c)=b\circ (a\circ c), \; \forall a,b,c\in A.$$

\noindent Such algebras are called {\it Novikov algebras} \cite{Nov}.
Notice that Novikov algebra $W_1(m)$ is unital.
If right-symmetric algebra $A$ is Novikov algebra, we will use denotion
$A^{nov}.$

\subsection{Right-symmetric modules and comodules}
A vector space $M$ is said to be {\it module} over
right-symmetric algebra $A,$ if it is endowed by  right action
$$M\times A\rightarrow M,\quad (m,a)\mapsto m\circ a$$

\noindent and left action
$$A\times M\rightarrow M, \quad (a,m)\mapsto a\circ m,$$

\noindent such that
$$m\circ [a,b]-(m\circ a)\circ b+(m\circ b)\circ a=0,$$
$$(a\circ m)\circ b-a\circ(m\circ b)-(a\circ b)\circ m+a\circ(b\circ m)=0,$$

\noindent for any $a,b\in A, m\in M.$ We will say, that $M$ is {\it antisymmetric}
$A-$module, if the left action of $A$ is trivial, i.e.,  $a\circ m=0,$ for
any $a\in A, m\in M.$ For module $M$ over right-symetric algebra $A,$ denote
by $M_{anti}$ its antisymmetric $A-$module: $M_{anti}=M,\; (m,a)\mapsto m\circ a, \;
(a,m)\mapsto 0,$ for all $m\in M_{anti}, a\in A.$

A right-symmetric $A-$module $M$ is said to be {\it special,}
if the right action  satisfies the following condition
$$m\circ (a\circ b)-(m\circ a)\circ b=0,\;\; \forall a,b\in A,\forall m\in M.$$
A special module is {\it antisymmetric,} if $a\circ m=0,$ for all $a\in A.$

{\bf Example.} For right-symmetric algebra $A$ its vector space $A$ can be
endowed by a natural structure of $A-$module,
$(a,m)\mapsto a\circ m, \; (m,a)\mapsto m\circ a,\; a,m\in A.$ In such cases
we say that $M=A$ is {\it regular} $A-$module. If $A$ is associative algebra,
then regular module is special.

The functor
$$A^{lie}\mbox{-module} \rightarrow \mbox{\; right\;} A^{lie}\mbox{-module}\rightarrow
\mbox{\;Antisymmetric\;} A\mbox{-module}$$

\noindent gives us an equivalence of the  category of antisymmetric $A-$modules
to the category of (right) $A^{lie}-$modules. Antisymmetric
$A$~-module corresponding to right $A^{lie}$~-module $M$ will be denoted
by $M_{anti}.$

Assume that $A$ is an associative algebra $A$ with multiplication $(a,b)\mapsto a\cdot b.$
In the last case of $A^{rsym}$ right-symmetric multiplications can be defined in two ways:
by $(a,b)\mapsto a\cdot b$ or by $(a,b)\mapsto b\cdot a.$
For definiteness we endow $A^{rsym}$ by multiplication
$(a,b)\mapsto a\cdot b.$
For associative algebra $A$ the functor
$$\mbox{Right\;} A^{ass}\mbox{-module}\rightarrow
\mbox{\;Antisymmetric special\;} A^{rsym}\mbox{-module}$$

\noindent gives us an equivalence of the categories
of antisymmetric special $A^{rsym}$~-modules and right $A^{ass}$-modules.

Right-symmetric $A$~-module $M$ can be endowed by a structure of module
over Lie algebra $A^{lie}$ by action $[a,m]=a\circ m-m\circ a.$  The 
obtained module is denoted by $M^{lie}.$

So, for defining module structure on a vector space $M$ over a
right-symmetric algebra $A$ one should define on $M$ right module structure
over the Lie algebra $A^{lie}$ and endow it by a left action that satisfies
condition $(AAM).$ As we mentioned before the last can be done by a trivial way
by setting $a\circ m=0, \forall a\in A,\forall m\in M.$

For a  module $M$ over right-symmetric algebra $A$ the subspace
$$M^{l.ass}=\{m:\in M : (m,a,b)=0,\forall a,b\in A\}$$

\noindent is called a {\it left associative invariant} subspace of $M,$ and
$$M^{l.inv}=\{m\in M: m\circ a=0, \forall a\in A\}$$

\noindent is called a {\it left invariant} subspace of $M.$ If $M=A$ is
{\it regular module,} then $A^{l.ass}$ is called a {\it left associative
center.} Notice that, $A^{l.inv}$ coincides with the left center of $A.$
Notice that, $M^{l.inv}$ is close under right action of $A$ and
$$M^{l.inv}\subseteq M^{l.ass}.$$

Module $M$ of (left) unital right-symmetric algebra is called
{\it (left) unital,} if
$$e\circ m=m,\quad \forall e\in Q_l(A), \forall m\in M.$$

\noindent  and {\it (left) central,} if
$$z\circ m=0, \quad, \forall z\in Z_{l}(A), \forall m\in M.$$

Regular module of unital right-symmetric algebra is unital and central.

A vector space $M$ is called {\it comodule} over
right-symmetric algebra $A,$ if there are given right action
$$M\times A\rightarrow M,\quad (m,a)\mapsto m\circ a,$$

\noindent and left action
$$A\times M\rightarrow M, \quad (a,m)\mapsto a\circ m,$$

\noindent such that
$$[a,b]\circ m-a\circ (b\circ m)+b\circ (a\circ m)=0,$$
$$-b\circ(m\circ a)+(b\circ m)\circ a-m\circ(a\circ b)+(m\circ a)\circ b=0,$$

\noindent for any $a,b\in A, m\in M.$
A comodule $M$ is {\it special,} if it satisfies the identity
$$(a\circ b)\circ m=a\circ (b\circ m), \;\; \forall a,b\in A,\forall m\in M.$$
A (special) comodule $M$ is called {\it antisymmetric,}
if $m\circ a=0,$ for any $a\in A.$

{\bf Example.} Let $A$ be right-symmetric algebra, $M$ be $A-$module and
$M'=\{f:M\rightarrow {\cal K}\}$ be a space of linear functions on $M.$
Set
$$(a\circ f)(m)=f(m\circ a),\; (f\circ a)(m)=f(a\circ m).$$

\noindent Then $M'$ under actions $(a,f)\mapsto a\circ f,\;
(f,a)\mapsto f\circ a,$ can be endowed by a structure of $A-$comodule.
Check it.
$$\{[a,b]f-a\circ (b\circ f)+b\circ(a\circ f)\}(m)=$$
$$f(m\circ[a,b]-(m\circ a)\circ b+(m\circ b)\circ a)=0,$$

$$\{-b\circ(f\circ a)+(b\circ f)\circ a
-f\circ(a\circ b)+(f\circ a)\circ b\}(m)=$$
$$f(-a\circ(m\circ b)+ (a\circ m)\circ b
-(a\circ b)\circ m+a\circ(b\circ m))=0.$$

\noindent The $A-$comodule $A'$ for regular module $A$ is called
{\it coregular} comodule of $A.$ If $A$ is associative, then $A'$
is special comodule.

For $A-$comodule $M$ let
$$M^{r.ass}=\{m\in M: (a,b,m)=0,\;\forall a,b\in A\}$$

\noindent be a {\it right associative invariant} subspace of $M$ and
$$M^{r.inv}=\{m\in M: a\circ m=0,\;\forall a\in A\}$$

\noindent be a {\it right invariant} subspace of $M.$ Notice that
$M^{r.inv}$ is close under left action of $A.$ An inclusion takes place 
$$M^{r.inv}\subseteq M^{r.ass}.$$

\subsection{Antisymmetric module $C^1_{right}(A,M)$}\label{1-module}
\begin{prp} The space of linear maps $C^1_{rsym}(A,M):=C^1(A,M)=
\{f:A\rightarrow M\}$ can be endowed by a
structure of antisymmetric $A-$module, where the right action is given by
$$(f\circ a)(b)=f(b)\circ a-f(a\circ b)+b\circ f(a), \; a,b\in A.$$
\end{prp}

{\bf Proof.} For $f\in C^1(A,M), a,b,c\in A,$ we have
$$(f\circ [b,c])(a)-((f\circ b)\circ c)(a)+((f\circ c)\circ b)(a)=$$

$$d_{rsym}f(a,[b,c])-d_{rsym}([f,b])(a,c)+d_{rsym}([f,c])(a,b)=$$

$$a\circ f([b,c])-f(a\circ [b,c])+f(a)\circ [b,c]-$$
$$-a\circ [f,b](c)+[f,b](a\circ c)-[f,b](a)\circ c+$$
$$+a\circ [f,c](b)-[f,c](a\circ b)+[f,c](a)\circ b=$$

$$a\circ f([b,c])-f(a\circ [b,c])+f(a)\circ [b,c]$$
$$-a\circ d_{rsym}f(c,b)+d_{rsym}f(a\circ c,b)-d_{rsym}f(a,b)\circ c+$$
$$+a\circ d_{rsym}f(b,c)-d_{rsym}f(a\circ b,c)+d_{rsym}f(a,c)\circ b=$$

$$\mathop{a\circ f([b,c])}\limline-\mathop{f(a\circ [b,c])}\limeq+
\mathop{f(a)\circ [b,c]}\limequiv$$

$$-\mathop{a\circ (c\circ f(b))}\limsim+\mathop{a\circ f(c\circ b)}\limline
-\mathop{a\circ (f(c)\circ b)}\limsimeq$$
$$+\mathop{(a\circ c)\circ f(b)}\limsim-\mathop{f((a\circ c)\circ b)}\limeq+
\mathop{f(a\circ c)\circ b}\limasymp$$
$$-\mathop{a\circ f(b)\circ c}\limsim+\mathop{f(a\circ b)\circ c}\limcong-
\mathop{(f(a)\circ b)\circ c}\limequiv+$$

$$+\mathop{a\circ (b\circ f(c))}\limsimeq-\mathop{a\circ f(b\circ c)}\limline+
\mathop{a\circ(f(b)\circ c)}\limsim$$
$$-\mathop{(a\circ b)\circ f(c)}\limsimeq+\mathop{f((a\circ b)\circ c)}\limeq-
\mathop{f(a\circ b)\circ c}\limcong+$$
$$\mathop{(a\circ f(c))\circ b}\limsimeq-\mathop{f(a\circ c)\circ b}\limasymp+
\mathop{(f(a)\circ c)\circ b}\limequiv=$$

$$=0.$$

\noindent So, $C^1(A,M)$ is a right $A^{lie}$-module. $\bullet$

\subsection{Universal enveloping algebras of right-symmetric algebras}
\label{universal}

Consider two copies of $A,$ denote them by $A^{r}, A^{l},$ and the tensor
algebra $T(A^{r}\oplus A^{l}).$  Algebras $A^{r}, A^{l}$ supposed to be
free as ${\cal K}-$module and the tensor algebra $T(A^{r}\oplus A^{l})$ is
associative and unital. Elements of $A^{r}$ and $A^{l}$ corresponding
to $a\in A$ denote as $r_a$ and $l_a.$ Let $U(A)$ be a factor-algebra of
$T(A^{r}\oplus A^{l}$ over an ideal $J$ generated by
$r_{[a,b]}-r_ar_b+r_br_a,\; [r_b,l_a]-l_bl_a+l_{a\circ b}.$
This algebra can be considered as {\it a universal enveloping algebra} of
right-symmetric algebra $A.$ Denote by $\breve{U}(A)$ the factor-algebra of
$T(A^{r}\oplus A^{l})$ over an ideal $\breve J$ generated by $\{r_{a\circ b}
-r_ar_b,\;[r_b,l_a]-l_bl_a+l_{a\circ b}\}.$ This algebra is called a {\it
special universal enveloping algebra of $A.$} Notice that,
$J\subset\breve J,$ since
$$r_{[a,b]}-[r_a,r_b]=\{r_{a\circ b}-r_{a}r_{b}\}\;-\;
\{ r_{b\circ a}+r_br_a\}\in \breve J.$$

\noindent So, the following exact sequences of algebras take place
$$0\rightarrow J\rightarrow T(A^{r}\oplus A^{l})\rightarrow U(A)\rightarrow 0,$$
$$0\rightarrow \breve J\rightarrow T(A^{r}\oplus A^{l})\rightarrow \breve U(A)
\rightarrow 0,$$

\noindent and
$$0\rightarrow \breve J/J\rightarrow U(A)\rightarrow \breve U(A)\rightarrow 0.$$

\noindent In particular, we can consider $\breve U(A)$ as right $U(A)-$module:
$$\breve u\bar v=\breve{u}\breve{v},$$

\noindent where $\breve u$ and $\bar u$ are elements of $\breve U(A)$ and $U(A)$
corresponding to $u\in T(A^{r}\oplus A^{l}).$

\begin{th} Let $A$ be a right-symmetric algebra.

i) There exists an equivalence of the categories of $A-$modules
and right $U(A)-$modules. The same is true for $A-$comodules
and left $U(A)-$modules.

ii) The category of special $A-$modules is equivalent to the category of
right $\breve{U}(A)-$modules. The same is true for special $A-$comodules
and left $\breve{U}(A)-$modules.
\end{th}

{Proof.} i) Let $(r,l): A\rightarrow End\, M$ be a representation of
right-symmetric algebra $A$ corresponding to $A$~-module $M,$ i.e.,
$$r:A\rightarrow End M, \;\; a\mapsto r_a, \;\; mr_a=m\circ a,$$
$$l:A\rightarrow End M, \;\; a\mapsto l_a, \;\; ml_a=a\circ m,$$

\noindent linear operators, such that for any $a,b\in A,$
$$r_{[a,b]}-r_ar_b+r_br_a=0, \eqno{(MAA)}$$
$$[r_b,l_a]-l_bl_a+l_{a\circ b}=0. \eqno{(AAM)}$$

\noindent So, any $A-$module is a right $U(A)-$module and, converse,
any right $U(A)-$module can be considered as an $A-$module.

A corepresentation $(r^{co},l^{co}):A\rightarrow End\, M,$
corresponding to $A\mbox{-comodule\;}$ $M,$
$$r^{co}:A\rightarrow End\, M, \; a\mapsto r_a, \; r^{co}_am=a\circ m,$$
$$l^{co}:A\rightarrow End\, M, \; a\mapsto l_a, \; l^{co}_am=m\circ a,$$

\noindent satisfies conditions (MAA), (AAM) for $r^{co}_a, l^{co}_a.$ So,
any $A-$comodule is a left $U(A)-$module. Any left $U(A)-$module can be
considered as a $A-$comodule.

ii) Let $M$ be a special $A-$module. Then by the rule
$$ mr_a=m\circ a, ml_a=m\circ a,$$

\noindent we obtain a right $\breve{U}(A)-$module:
$$m\circ(a\circ b)-(m\circ a)\circ b=0 \rightarrow
r_{a\circ b}=r_ar_b.$$

\noindent Converse, for a right $\breve{U}(A)-$module $N$,
one can correspond special $A-$module $N,$ by $n\circ a:=nr_a,
a\circ n=nl_a.$

For a special $A-$comodule $M$ notice that
$$(a\circ b)\circ m-a\circ (b\circ m)=0\Rightarrow
r^{co}_{a\circ b}=r^{co}_ar^{co}_b,$$

\noindent if $r^{co}_am=a\circ m, l^{co}_am=m\circ a.$ So, any
special $A-$comodule is a left $U^{spec}(A)-$module.
A converse statement is also evident. $\bullet$

\subsection{Right-symmetric cohomologies as a derived functor \label{derived}}

Recall that factor-images of the element $u\in T(A^{r}\oplus A^{l})$
in $U(A)$ and $\breve U(A)$ are denoted by $\bar u$ and $\breve u.$
Consider  $\breve A=A\oplus <1>$  as a right $U(A)-$module:
$$1\circ \bar r_a=r_a\circ 1=a,\;
1\circ \bar l_a=l_a\circ 1=a,\; a\circ \bar r_b=a\circ b, \;
a\circ \bar l_b=b\circ a,$$

\noindent for all $a\in A.$ Endow $\breve U(A)$ by a structure of
$U(A)-$module as in the subsection~\ref{universal}.
Consider $A\otimes\wedge^k(A)\otimes U(A), k\ge 0,$ as a right $U(A)-$module.
Then $A\otimes \wedge^kA\otimes U(A)$ is a free $U(A)-$module. Denote its
generators by $<a_0,a_1,\ldots,a_k>,$ where $a_0\in A,a_1\wedge\cdots\wedge
a_k\in \wedge^kA.$
Construct homomorphisms
$$\der:A\otimes \wedge^kA\otimes U(A)\rightarrow A\otimes\wedge^{k-1}A\otimes U(A), k>0,$$
$$\der:A\otimes U(A)\rightarrow \breve U(A)$$
$$\epsilon:\breve U(a)\rightarrow \breve A,$$

\noindent as below
$$\der<a_0,a_1,\ldots,a_k>=$$
$$\sum_{i=1}^k\{(-1)^{i+1}<a_i,a_1,\ldots,\hat{a_i},\ldots,a_k>\bar l_{a_0}$$
$$+(-1)^i<a_0\circ a_i,a_1,\ldots,\hat{a_i},\ldots,a_k>$$
$$+(-1)^{i+1}<a_0,a_1,\ldots,\hat{a_i},\ldots,a_k>\bar r_{a_i}\}$$
$$+\sum_{i<j}(-1)^{i+1}<a_0,a_1,\ldots,\hat{a_i},\ldots,a_{j-1},
[a_i,a_j],\ldots,a_k>,$$

$$\der(a_0)=<\breve 1>\bar l_{a_0}-<\breve 1>\bar r_{a_0},$$

$$\epsilon{\breve 1}=\breve 1, \;\;\epsilon{\breve r_a}=a,
\;\;\epsilon{\breve l_a}=-a.$$

Then the following sequence
$$\cdots\stackrel{\der}{\rightarrow} A\otimes\wedge^2A\otimes U(A)
\stackrel{\der}{\rightarrow}A\otimes A\otimes U(A)
\stackrel{\der}{\rightarrow} A\otimes U(A)\stackrel{\der}{\rightarrow} \breve{U}(A)
\stackrel{\epsilon}{\rightarrow} \breve A\rightarrow 0$$

\noindent is almost a free resolution of the right $U(A)-$module $\breve A.$
Here the words "almost free" mean that all members of the resolution except
$\breve U(A)$ is are free right $U(A)-$modules.

Notice that
$$Hom_{U(A)}(A\otimes\wedge^kA\otimes U(A),M)\cong A\otimes \wedge^kA,
\;k\ge 0,$$

\noindent and $Hom_{U(A)}(\breve U(A),M)$ consists of 
$g:\breve U(A)\rightarrow M,$ such that
$$g(\breve 1)(\bar r_a\bar r_b-\bar r_{a\circ b})=
g(\breve 1(\bar r_a\bar r_b-\bar r_{a\circ b}))=$$
$$g(\breve {r_ar_b}-\breve{r_{a\circ b}})=
g(\breve {r_a}{\breve r_b}-\breve{r_{a\circ b}})=$$
$$0,$$

\noindent for any $a,b\in A.$ So,
$$Hom_{U(A)}(\breve U(A),M)\cong \{m\in M: (m,a,b)=0\}.$$

Therefore, as a right-symmetric cochain complex we can take
$$C^*_{rsym}(A,M)=\oplus_kC^k_{rsym}(A,M),$$
$$C^0_{rsym}(A,M)=\{m\in M: (m,a,b)=0,\forall a,b\in A\},$$
$$C^{k+1}_{rsym}(A,M)=A\otimes \wedge^kA,\;\; k\ge 0.$$

\noindent These statements will follow from our results on right-symmetric
cohomologies in the next sections. Our approach is slightly different from Koszul's
approach. We will argue in cohomological terms and prove that right-symmetric
cochain complex has a pre-simplicial structure.

Let us mention these results relating homologies.
Let $M$ be comodule over right-symmetric algebra $A.$ Endow $M\otimes A$ by
a structure of antisymmetric $A-$comodule with a left action
$$b\circ(m\otimes a)=m\circ a\otimes b-m\otimes a\circ b+b\circ m\otimes a.$$

\noindent Set
$$C_0^{rsym}(A,M):=M^{r.ass}:=\{m\in M: (a,b,m)=0,\forall a,b\in A\},$$
$$C_{k+1}^{rsym}(A,M)=M\otimes A\otimes \wedge^k(A),\;\; k\ge 0.$$
$$C_*^{rsym}(A,M)=\oplus_kC^{rsym}_k(A,M).$$

\noindent Then $C_*^{rsym}(A,M)$ is chain complex under the boundary operator

$$\der: C_{k+1}^{rsym}(A,M)\rightarrow C_{k}^{rsym}(A,M),$$

$$\der(m\otimes a_0\otimes a_1\wedge\cdots\wedge a_k)=$$

$$\sum_{i=1}^{k}
(-1)^{i+1}\{m\circ a_0\otimes a_i\otimes a_1\wedge \cdots\hat{a_i}
\cdots\wedge a_k$$
$$-m\otimes a_0\circ a_i\otimes a_1\wedge \hat{a_i}\cdots\wedge a_k$$
$$+a_i\circ m\otimes a_0\otimes a_1\wedge\cdots\hat{a_i}\cdots\wedge a_k\}$$
$$+\sum_{i<j}(-1)^{i+1}m\otimes a_0\otimes a_1\wedge\cdots\hat{a_i}\cdots\wedge
a_{j-1}\wedge[a_i,a_j]\wedge\cdots a_k.$$

\noindent Moreover, $C_*^{rsym}(A,M)$ has antisymmetric $A-$comodule structure
with left action

$$\rho^{rsym}_{co}(x):C_{k+1}^{rsym}(A,M)\rightarrow C^{rsym}_{k+1}(A,M),$$

$$\rho^{rsym}_{co}(x)(m\otimes a_0\otimes a_1\wedge\cdots\wedge a_k)=$$

$$\sum_{i=1}^k(-1)^{i}(m\circ a_0\otimes a_i\otimes a_1\wedge\cdots\hat{a_i}\cdots
\wedge a_k$$
$$-m\otimes a_0\circ a_i\otimes a_1\wedge\cdots\hat{a_i}\cdots\wedge a_k$$
$$+a_i\circ m\otimes a_0\otimes a_1\wedge\cdots\hat{a_i}\cdots\wedge a_k\}$$
$$+\sum_{i<j}(-1)^{i+1}m\otimes a_0\otimes a_1\wedge \cdots\wedge a_{i-1}
\wedge [x,a_i]\wedge\cdots\wedge a_k,$$

\noindent and an isomorphism of $A-$comodules takes place
$$C_{k+1}^{rsym}(A,M)\cong C_k^{lie}(A,M\otimes A).$$

\noindent that induces an isomorphism of homology spaces
$$H_{k+1}^{rsym}(A,M)\cong H_k^{lie}(A,M\otimes A), \;\; k>0.$$

\subsection{Comultiplication of universal enveloping algebra}\label{comult}
Let $U(A)$ be the universal enveloping algebra of a right-symmetric algebra $A.$
As we noticed in section~\ref{1}, it can be generated by the elements $r_a, l_a, a\in A,$
such that
$$r_{[a,b]}-[r_a,r_b]=0,\; [l_a,r_b]-l_{a\circ b} +l_bl_a=0, \; a,b\in A.$$

\noindent Define homomorphism
$$\Delta:U(A)\rightarrow U(A)\otimes U(A),$$

\noindent by
$$\Delta(1)=1\otimes 1,$$
$$\Delta(r_a)=r_a\otimes 1+1\otimes (r_a-l_a),$$
$$\Delta(l_a)=l_a\otimes 1.$$

\noindent Since, according to right-symmetry identities,
$$\Delta([r_a,r_b])=$$
$$\Delta(r_a)\Delta(r_b)-\Delta(r_b)\Delta(r_a)=$$
$$r_ar_b\otimes 1+1\otimes (r_a-l_a)(r_b-l_b)-r_br_a\otimes 1+1\otimes (r_b-l_b)(r_a-l_a)=$$
$$r_{[a,b]}\otimes 1+1\otimes (r_{[a,b]}-l_{[a,b]}=$$
$$\Delta(r_{[a,b]}),$$

$$\Delta([l_a,r_b]-l_{a\circ b}+l_bl_a)=$$
$$(l_a\otimes 1)(r_b\otimes 1+1\otimes (r_b-l_b))-
(r_b\otimes 1+1\otimes (r_b-l_b))(l_a\otimes 1)-l_{a\circ b}\otimes 1+l_bl_a\otimes 1=$$
$$([l_a,r_b]-l_{a\circ b}+l_bl_a)\otimes 1=$$
$$0.$$

\noindent this definition is correct.

\begin{th} For a right-symmetric algebra $A$ and its universal enveloping algebra $U(A)$
the following diagram is commutative
$$\begin{array}{ccc}
U(A)&\stackrel{\Delta}\rightarrow&U(A)\otimes U(A)\\
\downarrow\lefteqn{\Delta}&&\downarrow\lefteqn{\Delta\otimes 1}\\
U(A)\otimes U(A)&\stackrel{1\otimes \Delta}\rightarrow &U(A)\otimes U(A)\otimes U(A)\\
\end{array}$$
\end{th}

{\bf Proof.}We must check that
$$(1\otimes \Delta)\Delta(u)=(\Delta\otimes 1)\Delta (u),\; \forall u\in U(A).$$

\noindent We have
$$(1\otimes \Delta)\Delta(r_a)=$$
$$(r_a\otimes 1\otimes 1+1\otimes r_a\otimes 1+1\otimes 1\otimes (r_a-l_a))
-1\otimes l_a\otimes 1=$$
$$(r_a\otimes 1\otimes 1+1\otimes r_a\otimes 1 -1\otimes l_a\otimes 1)+
1\otimes 1\otimes (r_a-l_a)=$$
$$\Delta(r_a)\otimes 1+1\otimes (r_a-l_a)=$$
$$(\Delta\otimes 1)\otimes \Delta(r_a),$$

$$(1\otimes \Delta)\Delta(l_a)=l_a\otimes 1\otimes 1=(\Delta\otimes 1)
\Delta(l_a).\; \bullet$$

Similarly, homomorphism $\Delta_1$ defined below is also comultiplication,
$$\Delta_1:U(A)\rightarrow U(A)\otimes U(A),$$
$$\Delta_1(1)=1\otimes 1,$$
$$\Delta_1(r_a)=(r_a-l_a)\otimes 1+1\otimes r_a,$$
$$\Delta_1(l_a)=1\otimes l_a$$

So, we can construct for given $A-$modules $M$ and $N$ their
tensor products $M\otimes N$ with a module structure unduced by
comultiplication $\Delta:$
$$(m\otimes n)\circ a=m\circ a\otimes n+m\otimes [n,a],$$
$$a\circ (m\otimes n)=a\circ m\otimes n.$$

\noindent Moreover, it is possible for right-symmetric $A-$module $M$
and for $A^{lie}-$module $N.$ These module structures on tensor products
are associative: if $M,N,S$ are modules over right symmetric algebra $A,$
then
$$(M\otimes N)\otimes S\cong M\otimes (N\otimes S).$$

{\bf Definition.} {\it For given modules $M,N$ over right-symmetric algebra $A,$
a homomorphism of $A-$modules $M\otimes N\rightarrow S$ is called a cup product of
$M$ and $N.$}

Denote the image of $m\otimes n$ in $S$ by $m\cup n.$ Thus, a bilinear map
$$M\times N\rightarrow S,\; (m,n)\mapsto m\cup n,$$

\noindent is said to be the cup product (pairing) of $M$ and $N$ to $S,$ if
$$(m\cup n)\circ a=m\circ a\,\cup n + m\cup [n,a],$$
$$a\circ (m\cup n)=a\circ m\,\cup n,$$

\noindent for any $a\in A, m,n\in M.$

Let
$$C^1(A,M)=Hom_k(A,M),\; C^{k}_{lie}(A,M)=Hom_k(\wedge^k A,M),\;k\ge 0,$$
$$C^{k+1}_{rsym}(A,M)=Hom_k(A\otimes \wedge^k(A),M),\;k\ge 0.$$

\begin{prp}\label{rho} $C^{k+1}_{rsym}(A,M)$ has an antisymmetric $A-$module
structure, where the right action
$$(C^{k+1}_{rsym}(A,M)\times A\rightarrow C^{k+1}_{rsym}(A,M),\;\;
(\psi,x)\mapsto \psi\circ x,$$

\noindent is defined by
$$(\psi\circ x)(a_0,a_1,\ldots,a_{k})=$$
$$a_0\circ \psi(x,a_1,\ldots,a_{k})-\psi(a_0\circ x,a_1,\ldots,a_{k})$$
$$+\psi(a_0,a_1,\ldots,a_k)\circ x+\sum_{i=1}^k\psi(a_0,a_1,\ldots,a_{i-1},[x,a_i],
\ldots,a_{k}),$$

\noindent for $\psi\in C^{k+1}_{rsym}(A,M),\; k\ge 0.$
\end{prp}

{\bf Proof.} Since,
$$C^1_{rsym}(A,M)=C^1(A,M),$$

\noindent an isomorphism of linear spaces takes place
$$G:C^1_{rsym}(A,M)\otimes C^k_{lie}(A,k)\rightarrow C^{k+1}_{right}(A,M),\; k\ge 0,$$
$$(G(f\otimes\psi))(a_0,a_1,\ldots,a_k)=f(a_0)\psi(a_1,\ldots,a_k).$$

\noindent In section~\ref{1-module} we have constructed an antisymmetric
right-module structure on $C^1_{right}(A,M).$ Lie module structure on
$C^k_{lie}(A,k)$ over $A^{lie}$ is well known. So, for
an antisymmetric $A-$module structure
$$C^1_{rsym}(A,M)\otimes C^k_{lie}(A,k)=\{f\otimes\phi: f\in C^1_{right}(A,M),
\phi\in C^k_{lie}(A,k)\}$$

\noindent we have
$$((f\otimes \phi)\circ x)(a_0\otimes (a_1,\ldots,a_k))=$$

$$((f\circ x)\otimes\psi)(a_0\otimes (a_1,\ldots,a_k))+
(f\otimes [\psi,x])(a_0\otimes (a_1,\ldots,a_k))=$$

$$(f(a_0)\circ x-f(a_0\circ x)+a_0\circ f(x))\otimes \psi(a_1,\ldots,a_k)$$
$$+f(a_0)\otimes \sum_{i=1}^k\psi(a_1,\ldots, [x,a_i],\ldots,a_k)=$$

\noindent We see that
$$G\{(f\otimes\psi)\circ x\}=\{G(f\otimes\psi)\}\circ x.$$

\noindent Therefore, $(\psi,x)\mapsto \psi\circ x$ gives us a right
representation.$\bullet$

\subsection{Right-symmetric modules for $W_n^{rsym}$}
Let
$$\Gamma_n=\{\alpha=(\alpha_1,\ldots,\alpha_n), \alpha_i\in{\bf Z},
i=1,\ldots, n\}$$
$$\Gamma_n^+=\{\alpha\in \Gamma_n: \alpha_i\ge 0, i=1,\ldots,n\}.$$

\noindent and
$$\Gamma_n({\bf m})=\{\alpha\in \Gamma_n^+: \alpha_i<p^{m_i}, i=1,\ldots,n\},$$

\noindent if $p>0, {\bf m}=(m_1,\ldots, m_n).$

Let
$$U={\cal K}[[x^{\pm 1},\ldots, x_n^{\pm 1}]]=\{x^{\alpha}: \alpha\in \Gamma_n\},$$
$$U^+={\cal K}[[x^{1},\ldots, x_n^ 1]]=\{x^{\alpha}: \alpha\in \Gamma_n^+\},$$

\noindent if $p=0,$ and
$$U=O_n({\bf m})=\{x^{(\alpha)}: \alpha\in \Gamma_n({\bf m})\},$$

\noindent if $p>0.$

For $p=0,$ let
$A=W_n^{rsym},$ if $U={\cal K}[[x^{\pm 1},\ldots, x^{\pm 1}_n]],$ and
$A^+=W_n^{+rsym},$ if $U^+={\cal K}[[x_1,\ldots,x_n]].$
Let $A$ be $W_n^{rsym}({\bf m}),$ if $U=O_n({\bf m}),\; p>0.$
Algebras $A, A^+$ are right-symmetric and $U$ is associative commutative.

Notice that $U$ has
a structure of antisymmetric graded $A-$module. The right action is given by
$u\circ a\der_i=a\der_i(u).$  The gradings are given by
$$|x^{\alpha}|=\sum_i\alpha_i, \;\;\;
\alpha\in \Gamma_n \;\;(\mbox{or\;} \Gamma_n({\bf m}) \mbox{\; if\;} p>0),$$
$$U=\oplus_kU_k,\;\; U_k=\{u\in U: |u|=k\},$$
$$A=\oplus_{k} A_k, \;\; A_k=\{a\der_i: |a|=k+1, i=1,\ldots,n\},\;
A_k\circ A_l\subseteq A_{k+l},$$
$$ U\circ U_l\subseteq U_{k+l}, \;\;
U_k\circ A_l\subseteq U_{k+l}, k,l\in {\bf Z}.$$

\noindent Notice that $A_0\cong gl_n^{rsym}.$

Let ${\cal A}_0=\oplus_kA_k,$ and ${\cal A^+}_0=\oplus_k>0 A_k^+,$ if $p=0.$
Let $M$ be  ${\cal A}_0-$module, if $p>0,$ and ${\cal A}_0^+-$module
if $p=0.$ Define antisymmetric $A-$module structure on
$U\otimes M_0$  by~(see \cite{Dzhumavestnik})
$$(u\otimes m)\circ a\der_i=a\der_i(u)\otimes m+
\sum_{\beta\in \Gamma_n}u\der^{\beta}(a)\otimes [m,x^{(\beta)}\der_i], \;\; p>0,$$

$$(u\otimes m)\circ a\der_i=a\der_i(u)\otimes m+
\sum_{\beta\in \Gamma^+_n}(1/\beta !)u\der^{\beta}(a)\otimes [m,x^{\beta}\der_i], \; p=0.$$

\section{Cohomologies of right-symmetric algebras}

\subsection {Pre-simplicial structures on $C^{*+1}_{rsym}(A,M)$}

For a right-symmetric algebra $A$ and its module $M$ we introduce a
structure of pre-simplicial cochain complex on $C_{rsym}^{*+1}(A,M)=
\oplus_{k\ge 0} C^{k+1}_{rsym}(A,M),$ where
$$C^{k+1}_{rsym}(A,M)=Hom(A\otimes \wedge^k A, M), \; k\ge 0.$$

\noindent Define linear operators $D_i: C^{*+1}_{rsym}(A,M)\rightarrow
C^{*+1}_{rsym}(A,M), i=1,2,\ldots,$ by the rules
$$D_i:C^{k}_{rsym}(A,M)\rightarrow C^{k+1}_{rsym}(A,M),$$
$$D_i\psi(a_0,a_1,\ldots,a_k)=$$
$$a_0\circ\psi(a_i,a_1,\ldots,\hat{a_i},\ldots, a_k)-
\psi(a_0\circ a_i,a_1,\ldots,\hat{a_i},\ldots,a_k)$$
$$+\psi(a_0,a_1,\ldots,\hat{a_i},\ldots,a_k)\circ a_i+
\sum_{i<j}\psi(a_0,a_1,\ldots,\hat{a_i},\ldots, a_{j-1},{[a_i,a_j]},\ldots,a_k),$$
$0\le k,\; i\le k,$
$$D_i\psi=0, \; i>k.$$

\noindent Here $\hat{a}$ means that the element $a$ is omitted.

In the next section we will endow $C^*_{rsym}(A,M)=\oplus_kC^k_{rsym}(A,M)$
by a structure of cochain complex, where
$$C_{rsym}^k(A,M)=0, \quad k<0,$$
$$C^0_{rsym}(A,M) =\{m\in M: (ma)b=m(ab), \quad \forall a,b\in A\}.$$

\begin{th} The set of endomorphisms ${D_i,i=1,2,\ldots}$ endows
$C^{*+1}_{rsym}(A,M)=\oplus_{k>0}C^k_{rsym}(A,M)$ by a pre-simplicial structure:
$$D_jD_i=D_iD_{j-1}, \; i<j,$$

\noindent In particular, $d_{rsym}=-\sum_i(-1)^iD_i,$ is a coboundary operator on
$C_{rsym}^{*+1}(A,M):$
$${d_{rsym}}^2=0.$$
\end{th}

{\bf Proof.} For $i<j, 1<k,$ we have
$$D_jD_i\psi(a_0,a_1,\ldots,a_k)=X_1+X_2+X_3+X_4,$$

\noindent where
$$X_1=a_0(D_i\psi(a_j,a_1,\ldots,\hat{a_j},\ldots,a_k),$$
$$X_2=-D_i\psi(a_0\circ a_j,a_1,\ldots,\hat{a_j},\ldots,a_k),$$
$$X_3=D_i\psi(a_0,a_1,\ldots,\hat{a_j},\ldots,a_k)a_j,$$
$$X_4=\sum_{j<s}D_i\psi(a_0,\ldots,\hat{a_j},\ldots,a_{s-1},{[a_j,a_s]},\ldots,a_k).$$

\noindent Direct calculations show that
$$X_1=$$
$${\mathop{a_0(a_j\psi(a_i,a_1,\ldots,\hat{a_i},\ldots,\hat{a_j},\ldots,a_k))}\limline}$$
$${\mathop{-a_0\psi(a_j\circ a_i,a_1,\ldots,\hat{a_i},\ldots,\hat{a_j},
\ldots,a_k)}\limeq}$$
$${\mathop{+a_0(\psi(a_j,a_1,\ldots,\hat{a_i},\ldots,\hat{a_j},
\ldots,a_k)a_i)}\limequiv}$$
$${+\sum_{i<s, s\ne j}\mathop{a_0(\psi(a_j,a_1,\ldots,\hat{a_i},\ldots,[a_i,a_s],
\ldots,a_k))}\limsimeq},$$

$$X_2=$$
$${\mathop{-(a_0\circ a_j)\psi(a_i,a_1,\ldots,\hat{a_i},\ldots,\hat{a_j},\ldots,a_k)}
\limsim}$$
$${\mathop{+\psi((a_0\circ a_j)\circ a_i,a_1,\ldots,\hat{a_i},\ldots,\hat{a_j},\ldots,a_k)}
\limapprox}$$
$${\mathop{-(\psi(a_0\circ a_j,a_1,\ldots,\hat{a_i},\ldots,\hat{a_j},\ldots,a_k))a_i}\limcong}$$
$${-\sum_{i<s, s\ne j}\mathop{\psi(a_0\circ a_j,a_1,\ldots,\hat{a_i},\ldots,[a_i,a_s],
\ldots,a_k)}\limasymp},$$

$$X_3=$$
$${\mathop{+(a_0\psi(a_i,a_1,\ldots,\hat{a_i},\ldots,\hat{a_j},\ldots,a_k))a_j}\limline}$$
$${\mathop{-(\psi(a_0\circ a_i,a_1,\ldots,\hat{a_i},\ldots,\hat{a_j},
\ldots,a_k))a_j}\limsmile}$$
$${\mathop{+((\psi(a_0,a_1,\ldots,\hat{a_i},\ldots,\hat{a_j},\ldots,a_k)a_i)a_j}\limfrown}$$
$$ {+\sum_{i<s, s\ne j}\mathop{(\psi(a_0,a_1,\ldots,\hat{a_i},\ldots,
a_{s-1},[a_i,a_s],\ldots,a_k))a_j}\limdoteq}$$
$${+\sum_{j<s}\mathop{a_0(\psi(a_i,a_1,\ldots,\hat{a_i},\ldots,\hat{a_j},\ldots,a_{s-1},
[a_j,a_s],\ldots,a_k))}\limequi},$$

$$X_4=$$
$${-\sum_{j<s}\mathop{\psi(a_0\circ a_i,a_1,\ldots,\hat{a_i},\ldots,\hat{a_j},\ldots,a_{s-1},
[a_j,a_s],\ldots,a_k)}\limcon}$$
$${+\sum_{j<s}\mathop{(\psi(a_0,a_1,\ldots,\hat{a_i},\ldots,\hat{a_j},\ldots,a_{s-1},
[a_j,a_s],\ldots,a_k))a_i}\limsi}$$
$${+\sum_{j<s, i<s_1, s<s_1}
\mathop{\psi(a_0,a_1,\ldots,\hat{a_i},\ldots,\hat{a_j},\ldots,a_{s-1},
[a_j,a_s],\ldots,a_{s_1-1},[a_i,a_{s_1}],\ldots,a_k)}\limsime}$$
$${+\sum_{j<s, i<s_1, s_1< s, s_1\ne j}
\mathop{\psi(a_0,a_1,\ldots,\hat{a_i},\ldots,\ldots,a_{s_1-1},
[a_j,a_{s_1}],\ldots,a_{s-1},[a_i,a_{s}],\ldots,a_k)}\limli}$$
$${+\sum_{j<s, (i<s_1, s_1=s)}\mathop
{\psi(a_0,a_1,\ldots,\hat{a_i},\ldots,\hat{a_j},\ldots,a_{s-1},
[a_i,[a_j,a_s]],\ldots,a_k)}\lime}.$$

Analogously,
$$D_iD_{j-1}\psi(a_0,a_1,\ldots,a_k)=Y_1+Y_2+Y_3+Y_4,$$

\noindent where
$$Y_1=a_0(D_{j-1}\psi(a_i,a_1,\ldots,\hat{a_i},\ldots,a_k),$$
$$Y_2=-D_{j-1}(a_0\circ a_i,a_1,\ldots,\hat{a_i},\ldots,a_k),$$
$$Y_3=D_{j-1}\psi(a_0,a_1,\ldots,\hat{a_i},\ldots,a_k)a_i,$$
$$Y_4=\sum_{i<s}D_{j-1}\psi(a_0,a_1,\ldots,\hat{a_i},\ldots,a_{s-1},[a_i,a_s],
\ldots,a_k).$$

\noindent We have
$$Y_1=$$
$${\mathop{a_0(a_i\psi(a_j,a_1,\ldots,\hat{a_i},\ldots,\hat{a_j},\ldots,a_k))}
\limequiv}$$
$${\mathop{-a_0\psi(a_i\circ a_j,a_1,\ldots,\hat{a_i},\ldots,\hat{a_j},
\ldots,a_k)}\limeq}$$
$${\mathop{+a_0(\psi(a_i,a_1,\ldots,\hat{a_i},\ldots,\hat{a_j},
\ldots,a_k)a_j)}\limline}$$
$${+\sum_{j<s}\mathop{a_0\psi(a_i,a_1,\ldots,\hat{a_i},\ldots,\hat{a_j},\ldots,
a_{s-1},[a_j,a_s],\ldots,a_k)}\limequi},$$

$$Y_2=$$
$${\mathop{-(a_0\circ a_i)\psi(a_j,a_1,\ldots,\hat{a_i},\ldots,\hat{a_j},\ldots,a_k)}
\limsim}$$
$${\mathop{+\psi((a_0\circ a_i)\circ a_j,a_1,\ldots,\hat{a_i},
\ldots,\hat{a_j},\ldots,a_k)}\limapprox}$$
$${\mathop{-\psi(a_0\circ a_i,a_1,\ldots,\hat{a_i},\ldots,\hat{a_j},\ldots,a_k)a_j}
\limsmile}$$
$${-\sum_{j<s}\mathop{\psi(a_0\circ a_i,a_1,\ldots,\hat{a_i},
\ldots,\hat{a_j},\ldots,a_{s-1},[a_j,a_s],\ldots,a_k)}\limcon},$$

$$Y_3=$$
$${\mathop{(a_0\psi(a_j,a_1,\ldots,\hat{a_i},\ldots,\hat{a_j},\ldots,a_k))a_i}
\limequiv}$$
$${-\mathop{(\psi(a_0\circ a_j,a_1,\ldots,\hat{a_i},\ldots,
\hat{a_j},\ldots,a_k))a_i}\limcong}$$
$${+\mathop{((\psi(a_0,a_1,\ldots,\hat{a_i},\ldots,\hat{a_j},\ldots,a_k)a_j)a_i}
\limfrown}$$
$${+\sum_{j<s}\mathop{(\psi(a_0,a_1,\ldots,\hat{a_i},\ldots,
\hat{a_j},\ldots, a_{s-1},[a_j,a_s],\ldots,a_k))a_i}\limsi}.$$

\noindent Present $Y_4$ as a sum
$$Y_4=Y_{4,1}+Y_{4,2}+Y_{4,3},$$

\noindent where
$$Y_{4,1}=\sum_{i<s<j}D_{j-1}\psi(a_0,a_1,\ldots,\hat{a_i},\ldots,a_{s-1},[a_i,a_s],
\ldots,a_k),$$
$$Y_{4,2}=
D_{j-1}\psi(a_0,a_1,\ldots,\hat{a_i},\ldots,a_{j-1},[a_i,a_j],
\ldots,a_k),$$
$$Y_{4,3}=\sum_{j<s}D_{j-1}\psi(a_0,a_1,\ldots,\hat{a_i},\ldots,a_{s-1},[a_i,a_s],
\ldots,a_k).$$

\noindent These elements can be expressed in the following ways

$$Y_{4,1}=$$
$$\sum_{i<s<j}\mathop{a_0(\psi(a_j,a_1,\ldots,\hat{a_i},\ldots,a_{s-1},
[a_i,a_s],\ldots,\hat{a_j},\ldots, a_k))}\limsimeq$$
$$-\sum_{i<s<j}\mathop{\psi(a_0\circ a_j,a_1,\ldots,\hat{a_i},\ldots,a_{s-1},
[a_i,a_s],\ldots,\hat{a_j},\ldots,a_k)}\limasymp $$
$$+\sum_{i<s<j}\mathop{(\psi(a_0,a_1,\ldots,\hat{a_i},\ldots,a_{s-1},
[a_i,a_s],\ldots,\hat{a_j},\ldots,a_k))a_j}\limdoteq$$
$$+\sum_{i<s_1<j<s}\mathop{\psi(a_0,a_1,\ldots,\hat{a_i},\ldots,a_{s_1-1},
[a_i,a_{s_1}],\ldots,\hat{a_j},\ldots,a_{s-1},[a_j,a_{s}],\ldots,a_k)}
\limli,$$

$$Y_{4,2}=$$
$$\mathop{+a_0(\psi([a_i,a_j],a_1,\ldots,\hat{a_i},\ldots,
\hat{a_j},\ldots, a_k))}\limeq$$
$$-\mathop{\psi(a_0\circ [a_i,a_j],a_1,\ldots,\hat{a_i},
\ldots,\hat{a_j},\ldots,a_k)}\limapprox$$
$$+\mathop{(\psi(a_0,a_1,\ldots,\hat{a_i},\ldots,\hat{a_j},\ldots,
a_k))[a_i,a_j]}\limfrown$$
$$+\sum_{j<s}\mathop{\psi(a_0,a_1,\ldots,\hat{a_i},\ldots,
\hat{a_j},\ldots,a_{s-1},[[a_i,a_j],a_{s}],\ldots,a_k)}
\lime,$$

$$Y_{4,3}=$$
$$\sum_{j<s}\mathop{a_0(\psi(a_j,a_1,\ldots,\hat{a_i},\ldots,\hat{a_j},\ldots,a_{s-1},
[a_i,a_s],\ldots,a_k))}\limsimeq$$
$$-\sum_{j<s}\mathop{\psi(a_0\circ a_j,a_1,\ldots,\hat{a_i},\ldots,\hat{a_j},\ldots,a_{s-1},
[a_i,a_s],\ldots,a_k)}\limasymp$$
$$+\sum_{j<s}\mathop{(\psi(a_0,a_1,\ldots,\hat{a_i},\ldots,\hat{a_j},\ldots,a_{s-1},
[a_i,a_s],\ldots,a_k))a_j}\limeq$$
$$+\sum_{j<s<s_1}\mathop{\psi(a_0,a_1,\ldots,\hat{a_i},\ldots,\hat{a_j},\ldots,a_{s-1},
[a_i,a_s],\ldots,a_{s_1-1},[a_j,a_{s_1}],\ldots,a_k)}\limli$$
$$+\sum_{j<s}\mathop{\psi(a_0,a_1,\ldots,\hat{a_i},\ldots,\hat{a_j},\ldots,a_{s-1},
[a_j,[a_i,a_s]],\ldots,a_k)}\lime$$
$$+\sum_{j<s_1<s}\mathop{\psi(a_0,a_1,\ldots,\hat{a_i},\ldots,\hat{a_j},\ldots,a_{s_1-1},
[a_i,a_{s_1}],\ldots,a_{s-1},[a_j,a_{s}],\ldots,a_k)}\limsime$$

Using right-symmetric identity for the expressions underlined, in similar ways
we obtain that
$$D_jD_i=D_iD_{j-1},\qquad i<j. \qquad\qquad\bullet$$

\subsection{Cohomologies of right-symmetric algebras and Cartan's formulas.\label{cohom}}

In the  previous section we proved that $C^{*+1}_{rsym}(A,M)=\oplus_{k>0}C^k_{rsym}(A,M)$
is a cochain complex under coboundary operator $d_{rsym},$ such that

$$d_{rsym}\psi(a_0,a_1,\ldots,a_k)=$$
$$-\sum_{i=1}^{k}(-1)^{i}a_0\circ(\psi(a_i,a_1,\ldots,\hat{a_i},\ldots,a_{k}))$$
$$+\sum_{i=1}^{k}(-1)^{i}\psi(a_0\circ a_i,a_1,\ldots,\hat{a_i},\ldots,a_{k})$$
$$\sum_{1\le i<j\le k}(-1)^{i+1}\psi(a_0,a_1,\ldots,\hat{a_i},\ldots,[a_i,a_j],\ldots,a_{k})$$
$$-\sum_{i=1}^{k}(-1)^i(\psi(a_0,a_1,\ldots,\hat{a_i},\ldots,a_{k}))\circ a_i,$$
$$\psi\in C^{k}(A,M), 0<k.$$

\noindent For $m\in M,$ define $d_{rsym}\in C^1_{rsym}(A,M),$
$$d_{rsym}m(a)=a\circ m-m\circ a.$$

\noindent It is easy to see that
$$d_{rsym}^2m(a,b)=$$
$$a\circ d_{rsym}m(b)-d_{rsym}m(a\circ b)+d_{rsym}m(a)\circ b=$$
$$\mathop{a\circ (b\circ m)}\limline-
\mathop{a\circ(m\circ b)}\limeq-\mathop{(a\circ b)\circ m}\limline
+m\circ (a\circ b)+\mathop{(a\circ m)\circ b}\limeq-(m\circ a)\circ b=$$
$$(a,b,m)-(a,m,b)+m\circ(a\circ b)-(m\circ a)\circ b.$$

\noindent Thus, according to right-symmetric identity,
\begin{equation}
d_{rsym}^2m(a,b)=m\circ(a\circ b)-(m\circ a)\circ b. \label{cocycle}
\end{equation}

\noindent From this fact two conclusions follow.
Firstly, taking {\it a subspace of left associative invariants}
$$M^{l.ass}=\{m\in M: (m,a,b)=0, \; \forall a,b\in A\}$$

\noindent as a $0$-cochain subspace $C^0_{rsym}(A,M),$ we obtain
{\it cochain complex
$$C^*_{rsym}(A,M)=\oplus_{k\ge 0}C^k_{rsym}(A,M)$$

\noindent under coboundary operator $d_{rsym}.$} The second conclusion
will be discussed in the next section in the construction of standard 2-cocycles.

Let
$$Z^*_{rsym}(A,M)=\oplus_kZ^k_{rsym}(A,M),$$
$$ Z^k_{rsym}(A,M)=\{\psi\in C^k_{rsym}(A,M): d_{rsym}\psi=0\},$$

\noindent be spaces of right-symmetric cocycles,
$$B^*_{rsym}(A,M)=\oplus_kB^k_{rsym}(A,M),$$
$$ B^k_{rsym}(A,M)=\{d_{rsym}\omega: \omega\in C^{k-1}_{rsym}(A,M)\},$$

\noindent be spaces of right-symmetric coboundaries, and
$$H^*_{rsym}(A,M)=\oplus_kH^k_{rsym}(A,M),$$
$$H^k_{rsym}(A,M)=Z^k_{rsym}(A,M)/B^k_{rsym}(A,M),$$

\noindent be right-symmetry cohomology spaces.

{\bf Definitions.} { For any $x\in A$ the interior product endomorphism
$i(x)$ of $C^*_{rsym}(A,M)$ is defined by
$$i(x):C^{k+1}_{rsym}(A,M)\rightarrow C^{k}_{rsym}(A,M),$$
$$i(x)\psi(a_0,\ldots,a_{k-1})=\psi(a_0,x,a_1,\ldots,a_{k-1}),\quad k>0,$$
$$i(x)\psi=0,\quad \psi\in C^1_{rsym}(A,M).$$

\noindent Let $\rho_{lie}: A^{lie}\rightarrow C^{*+1}_{rsym}(A,M)$
be a representation of Lie algebra $A^{lie}$ corresponding to
antisymmetric representation
$$\rho_{rsym}:A\rightarrow C^{*+1}_{rsym}(A,M),
\; \rho_{rsym}(x)\psi=0,\; \psi\rho_{rsym}(x)=\psi\circ x,$$

\noindent constructed in proposition~\ref{rho},
$$(\rho_{lie}(x)\psi)(a_0,a_1,\ldots,a_{k})=$$
$$-a_0\circ \psi(x,a_1,\ldots,a_{k})+\psi(a_0\circ x,a_1,\ldots,a_{k})$$
$$-\psi(a_0,a_1,\ldots,a_k)\circ x-\sum_{i=1}^k\psi(a_0,a_1,\ldots,a_{i-1},
[x,a_i],\ldots,a_{k}).$$}

Recall that
$$\rho_{lie}(x)\psi=-\psi\rho_{lie}(x)=-\psi\circ x=-\psi\rho_{rsym}(x).$$

\begin{prp} {\rm (Cartan's formulas)}
Consider $C^{*+1}_{rsym}(A,M):=\oplus_kC^{k+1}_{rsym}(A,M)$
as a $A^{lie}-$module. For linear operators on $C^{*+1}_{rsym}(A,M)$
the following relations take place
$$i(x)D_l=D_{l-1}i(x), \; l>1, \leqno{(i)} $$
$$i(x)D_1=-\rho_{lie}(x), \leqno{(ii)}$$
$${\rho_{lie}}[x,y]=[{\rho_{lie}}(x),{\rho_{lie}}(y)], \leqno{(iii)}$$
$$[i(x),{\rho_{lie}}(y)]=-i([x,y]), \leqno{(iv)}$$
$$d_{rsym}i(x)+i(x)d_{rsym}=-{\rho_{lie}}(x), \leqno{(v)}$$
\end{prp}

{\bf Proof.} (i)$$i(x)D_l\psi(a_0,\ldots,a_{k-1})=$$

$$D_l\psi(a_0,x,a_1,\ldots,a_{k-1})=$$

$$a_0\circ\psi(a_{l-1},x,a_1,\ldots,\hat{a_{l-1}},\ldots, a_{k-1})$$
$$-\psi(a_0\circ a_{l-1},x,a_1,\ldots,\hat{a_{l-1}},\ldots,a_{k-1})$$
$$+\psi(a_0,x,a_1,\ldots,\hat{a_{l-1}},\ldots,a_{k-1})\circ a_{l-1}$$
$$+\sum_{l-1<j}\psi(a_0,x,a_1,\ldots,\hat{a_{l-1}},\ldots, a_{j-1},
{[a_{l-1},a_j]},\ldots,a_{k-1})=$$

$$a_0\circ i(x)\psi(a_{l-1},a_1,\ldots,\hat{a_{l-1}},\ldots, a_{k-1})$$
$$-i(x)\psi(a_0\circ a_{l-1},a_1,\ldots,\hat{a_{l-1}},\ldots,a_{k-1})$$
$$+i(x)\psi(a_0,a_1,\ldots,\hat{a_{l-1}},\ldots,a_{k-1})\circ a_{l-1}$$
$$+\sum_{l-1<j}i(x)\psi(a_0,a_1,\ldots,\hat{a_{l-1}},\ldots, a_{j-1},
{[a_i,a_j]},\ldots,a_{k-1})=$$

$$D_{l-1}i(x)\psi(a_0,\ldots,a_{k-1}).$$

\medskip

(ii) $$i(x)D_1\psi(a_0,\ldots,a_{k-1})=$$

$$D_1\psi(a_0,x,a_1,\ldots, a_{k-1})=$$

$$a_0\circ\psi(x,a_1,\ldots, a_{k-1})-
\psi(a_0\circ x,a_1,\ldots, a_{k-1})$$
$$+\psi(a_0,a_1,\ldots,a_{k-1})\circ x+
\sum_{0<j}\psi(a_0,a_1,\ldots, a_{j-1},{[x,a_j]},\ldots,a_{k-1})=$$

$$(\psi{\rho_{rsym}}(x))(a_0,\ldots,a_{k-1}).$$

\medskip

(iii) Proposition~\ref{rho}.

(iv) We obtain
$$-\{(i(x){\rho_{lie}}(y))\psi\}(a_0,a_1,\ldots,a_{k-1})=$$

$$ (\psi\circ y)(a_0,x,a_1,\ldots,a_{k-1})=$$

$$\mathop{a_0\circ(\psi(y,x,a_1,\ldots,a_{k-1})}\limline
-\mathop{\psi(a_0\circ y,x,a_1,\ldots,a_{k-1})}\limeq$$
$$+\mathop{\psi(a_0,x,a_1,\ldots,a_{k-1})\circ y}\limli$$
$$+\psi(a_0,[x,y],a_1,\ldots,a_{k-1})
+\sum_{i=1}^{k-1}\mathop{\psi(a_0,x,a_1,\ldots,
a_{i-1},[a_i,y],\ldots,a_{k-1})}\limsim,$$

\noindent and
$$-\{(\rho_{lie}(y)i(x))\psi\}(a_0,a_1,\ldots,a_{k-1})=$$

$$\{(i(x)\psi)\circ y\}(a_0,a_1,\ldots,a_{k-1})=$$

$$a_0\circ\{(i(x)\psi)(y,a_1,\ldots,a_{k-1})\}
-(i(x)\psi)(a_0\circ y,a_1,\ldots,a_{k-1})$$
$$+\{(i(x)\psi)(a_0,a_1,\ldots,a_{k-1})\}\circ y$$
$$+\sum_{i=1}^{k-1}(i(x)\psi)(a_0,a_1,\ldots,a_{i-1},[a_i,y],\ldots,a_{k-1})=$$

$$\mathop{a_0\circ(\psi(y,x,a_1,\ldots,a_{k-1}))}\limline
-\mathop{\psi(a_0\circ y,x,a_1,\ldots,a_{k-1})}\limeq$$
$$+\mathop{\psi(a_0,x,a_1,\ldots,a_{k-1})\circ y}\limli$$
$$+\sum_{i=1}^{k-1}\mathop{\psi(a_0,x,a_1,\ldots,a_{i-1},[a_i,y],
\ldots,a_{k-1})}\limsim.$$

\noindent Thus
$$\{(-i(x)\rho_{lie}(y)+\rho_{lie}(y)i(x))\psi\}(a_0,a_1,\ldots,a_{k-1})=$$

$$(i[x,y]\psi)(a_0,a_1,\ldots,a_{k-1}).$$

(v) According~(i) and~(ii),
$$i(x)d_{rsym}= i(x)D_1+\sum_{l>1}(-1)^{l+1}i(x)D_l=$$
$$-\rho_{lie}(x)+\sum_{l>1}(-1)^{l+1}D_{l-1}i(x)=$$
$$-\rho_{lie}(x)-\sum_{l>0}(-1)^{l+1}D_{l}i(x).$$

\noindent Thus,
$$i(x)d_{rsym}+d_{rsym}i(x)=-\rho_{lie}(x).\bullet $$

\subsection{Long exact cohomological sequence}\label{nabla}
The following theorem follows from standard homological results.

\begin{th} Let $A$ by a right-symmetric algebra and
$$0\rightarrow M\rightarrow T\rightarrow S\rightarrow 0$$

\noindent be a short exact sequence of right-symmetric $A$-modules. Then
an exact sequence of right-symmetric cohomology spaces take place
$$0\rightarrow Z^0_{rsym}(A,M)\rightarrow Z^0_{rsym}(A,T)\rightarrow
Z^0_{rsym}(A,S)\stackrel{\delta}\rightarrow $$
$$Z^1_{rsym}(A,M)\rightarrow Z^1_{rsym}(A,T)\rightarrow
Z^1_{rsym}(A,S)\stackrel{\delta}{\rightarrow} H^2_{rsym}(A,M)\rightarrow\cdots$$
$$\rightarrow H^k_{rsym}(A,M)\rightarrow H^k_{rsym}(A,T)\rightarrow H^k_{rsym}(A,S)
\stackrel{\delta}{\rightarrow}H^{k+1}_{rsym}(A,M)\rightarrow\cdots$$

\noindent Here $\delta$ is a connected homomorphism:
$$\delta[\psi]=[d_{rsym}\phi],\;\;  [\psi]\in H^{k}_{rsym}(A,S), k>1, $$
$$\delta \psi_1=[d_{rsym}\phi_1],\;\; \psi_1\in Z^i_{rsym}(A,S), i=0,1,$$

\noindent where $\phi\in Z^k_{rsym}(A,T)$  is a representative of
the cohomological class $[\psi]$ and $\phi_1\in Z^i_{rsym}(A,T)$ moves to
$\psi_1$ under a natural homomorphism $Z^i_{rsym}(A,T)\rightarrow
Z^i_{rsym}(A,S), i=0,1.$
\end{th}

Define a homomorphism $\nabla: S^{l.ass}\rightarrow Z^2_{rsym}(A,M),$
as a composition
$$\nabla:C^0_{rsym}(A,S)\stackrel{d_{rsym}}{\rightarrow} B^1_{rsym}(A,S)
\stackrel{\delta}{\rightarrow}
Z^2_{rsym}(A,M),\;\; s\mapsto d_{rsym}s\mapsto \delta (d_{rsym}s).$$

\noindent Then
$$\nabla({\tilde m})(a,b)={\tilde m}\circ (a\circ b)-({\tilde m}\circ a)\circ b.$$

\noindent In particular, there exist homomorphisms
$$\delta: S^{l.ass}\rightarrow Z^1_{rsym}(A,M),\;\; \delta(s):a\mapsto [a,s],$$
$$\delta: S^{l.inv}\rightarrow Z^1_{rsym}(A,M),\;\; -\delta(s):a\mapsto a\circ s.$$

\subsection{Connections between right-symmetric cohomologies and
Chevalley-Eilenberg cohomologies.}
Recall that for any $A^{lie}-$module $Q$ a standard representation
$\varrho: A^{lie}\rightarrow C^*_{lie}(A,Q)$ is given by
$$\varrho(x)\psi(a_1,\ldots,a_k)=$$
$$[x,\psi(a_1,\ldots,a_k)]-\sum_{i=1}^k\psi(a_1,\ldots,
a_{i-1},[x,a_i],\ldots,a_k),$$

\noindent where $(x,q)\mapsto [x,q]$ is representation corresponding to the 
Lie module $Q$ and $\psi\in C^k_{lie}(A,Q).$

\begin{th}\label{chev}
Let $A$ be a right-symmetric algebra and $M$ be an $A-$module.

\noindent An operator
\begin{equation}
F:C^{k}_{lie}(A,C^1(A,M))\rightarrow C^{k+1}_{rsym}(A,M), \quad k>0,\label{F}
\end{equation}

\noindent defined by the rule
$$F\psi(a_0,a_1,\ldots,a_k)=-\psi(a_1,\ldots,a_{k-1})(a_0)$$

\noindent induces an isomorphism of $A^{lie}-$modules. Moreover, $F$
induces an isomorphism of cochain complexes $C^{*+1}_{rsym}(A,M)$ and
$C^*_{lie}(A,C^1(A,M)).$ In particular,
\begin{equation}
H^{k+1}_{rsym}(A,M)\cong H^{k}_{lie}(A,C^1(A,M)), \; k>0. \label{a}
\end{equation}

\noindent The following sequence is exact
\begin{equation}
0\rightarrow Z^0_{rsym}(A,M)\rightarrow C^0_{rsym}(A,M)\rightarrow
H^0_{lie}(A,C^1(A,M))\rightarrow H^1_{rsym}(A,M)\rightarrow 0.\label{b}
\end{equation}
\end{th}

{\bf Proof.}
Prove that for any $x\in A, \; k>0,$ the following diagram is commutative

$$\begin{array}{ccc}
 C^k_{lie}(A,C^1(A,M))&\stackrel{\varrho(x)}\longrightarrow&C^{k}_{lie}(A,C^1(A,M))\\
\downarrow\lefteqn{F}&&\downarrow\lefteqn{F}\\
C^{k+1}_{rsym}(A,C^1(A,M))&\stackrel{\rho_{lie}(x)}\longrightarrow& C^{k+1}_{rsym}(A,M)\\
\end{array}$$

\medskip

\noindent For $\psi\in C^k_{lie}(A,C^1(A,M))$ we have
$$F\{\varrho(x)\psi\}(a_0,a_1,\ldots,a_{k+1})=
-\varrho(x)\psi(a_1,\ldots,a_{k+1})(a_0)=$$

$$-\{x\circ(\psi(a_1,\ldots,a_k))\}(a_0)+\sum_{i=1}^k\psi(a_1,\ldots,
a_{i-1},[x,a_i],\ldots,a_k)(a_0)=$$

$$+\{(\psi(a_1,\ldots,a_k))\circ x\}(a_0)+\sum_{i=1}^k\psi(a_1,\ldots,
a_{i-1},[x,a_i],\ldots,a_k)(a_0)=$$

$$-a_0\circ\psi(x,a_1,\ldots,a_k)+\psi(a_0\circ x,a_1,\ldots,a_k)-
\psi(a_0,a_1,\ldots,a_k)\circ x$$
$$-\sum_{i=1}^k\psi(a_0,a_1,\ldots,[x,a_i],\ldots,a_k)=$$

$$-(\psi\rho_{rsym}(x))(a_0,a_1,\ldots,a_k)=$$

$$\{(F\psi)\rho_{rsym}(x)\}(a_1,\ldots,a_k)(a_0).$$

\noindent Thus, $F:C^k_{lie}(A,C^1(A,M))\rightarrow C^{k+1}_{rsym}(A,M)$
is a homomorpism of $A^{lie}-$modules. It is evident that, $F$ has no kernel
and $F$ is epimorphism.

Now, prove that for any $k\ge 0$ the following diagram is commutative

$$\begin{array}{ccc}
 C^k_{lie}(A,C^1(A,M))&\stackrel{d_{lie}}\longrightarrow&C^{k+1}_{lie}(A,C^1(A,M))\\
\downarrow\lefteqn{F}&&\downarrow\lefteqn{F}\\
C^{k+1}_{rsym}(A,C^1(A,M))&\stackrel{d_{rsym}}\longrightarrow& C^{k+2}_{rsym}(A,M)\\
\end{array}$$

\medskip

\noindent For $\psi\in C^k_{lie}(A,C^1(A,M)),$
$$F(d_{lie}\psi)(a_0,a_1,\ldots,a_{k+1})=d_{lie}\psi(a_1,\ldots,a_{k+1})(a_0)=$$

$$\sum_{1\le i<j\le k+1}(-1)^{i}\psi(a_1,\ldots,\hat{a_i},\ldots,[a_i,a_j],\ldots,a_{k+1})(a_0)$$
$$-\sum_{i=1}^{k+1}(-1)^i[a_i,\psi(a_1,\ldots,\hat{a_i},\ldots,a_{k+1})](a_0)=$$

$$\sum_{1\le i<j\le k+1}(-1)^{i+1}F\psi(a_0,a_1,\ldots,\hat{a_i},\ldots,[a_i,a_j],\ldots,a_{k+1})$$
$$+\sum_{i=1}^{k+1}(-1)^id_{rsym}(\psi(a_1,\ldots,\hat{a_i},\ldots,a_{k+1}))(a_0,a_i)=$$

$$\sum_{1\le i<j\le k+1}(-1)^{i+1}F\psi(a_0,a_1,\ldots,\hat{a_i},\ldots,[a_i,a_j],\ldots,a_{k+1})$$
$$+\sum_{i=1}^{k+1}(-1)^ia_0\circ(\psi(a_1,\ldots,\hat{a_i},\ldots,a_{k+1}))(a_i)$$
$$-\sum_{i=1}^{k+1}(-1)^i(\psi(a_1,\ldots,\hat{a_i},\ldots,a_{k+1}))(a_0\circ a_i)$$
$$+\sum_{i=1}^{k+1}(-1)^i(\psi(a_1,\ldots,\hat{a_i},\ldots,a_{k+1})(a_0))\circ a_i=$$

$$\sum_{1\le i<j\le k+1}(-1)^{i+1}F\psi(a_0,a_1,\ldots,\hat{a_i},\ldots,[a_i,a_j],\ldots,a_{k+1})$$
$$-\sum_{i=1}^{k+1}(-1)^{i}a_0\circ(F\psi(a_i,a_1,\ldots,\hat{a_i},\ldots,a_{k+1}))$$
$$+\sum_{i=1}^{k+1}(-1)^{i}F\psi(a_0\circ a_i,a_1,\ldots,\hat{a_i},\ldots,a_{k+1})$$
$$-\sum_{i=1}^{k+1}(-1)^i(F\psi(a_0,a_1,\ldots,\hat{a_i},\ldots,a_{k+1}))\circ a_i$$

$$=d_{rsym}(F\psi)(a_0,a_1,\ldots,a_{k+1}).$$

\noindent Thus we obtain the equivalence
of cochain complexes $\oplus_{k>0}C^k_{lie}(A,C^1(A,M))$ and $\oplus_{k>1}C^k_{rsym}(A,M).$
In particular, an isomorphism~(\ref{a}) takes place.

Since, $H^1_{rsym}(A,M)=Z^1_{rsym}(A,M)/B^1_{rsym}(A,M),$ and
$$Z^1_{rsym}(A,M)=\{\psi\in C^1(A,M): d_{rsym}\psi=0\}=Z^0_{lie}(A,C^1(A,M)),$$
$$B^1_{rsym}(A,M)=\{d_{rsym}m:  m\in M, (m,a,b)=0,\forall a,b\in A\}=
M^{l.ass}/{M^{l.ass}\cap M^{inv}},$$

\noindent we have the exactness of (\ref{b}). $\bullet$

{\bf Definition.} Let $f:C^*_{rsym}(A,M)\rightarrow C^*_{lie}(A,M)$ be
a linear operator, such that
$$f:C^k_{rsym}(A,M)\rightarrow C^k_{lie}(A,M),$$
$$f\psi(a_1,\ldots,a_k)=\sum_{i=1}^k(-1)^{i+k+1}\psi(a_i,a_1,\ldots,\hat{a_i},\ldots,a_k).$$

\noindent Introduce subspaces
$$\bar{C}^*_{rsym}(A,M)=\oplus_k \bar{C}^k_{rsym}(A,M),$$

$$\bar{C}^k_{rsym}(A,M)= \{\psi\in C^{k+2}_{rsym}(A,M): f\psi=0\},\; k\ge 0 \; \bullet$$

\begin{th} Let $A$ be a right-symmetric algebra and $M$ be an $A$-module.
Then the operator $f:C^*_{rsym}(A,M)\rightarrow C^*_{lie}(A,M)$
is the homomophism of cochain complexes and the folllowing cohomological sequence

$$0\rightarrow Z^1_{rsym}(A,M)\rightarrow Z^1_{lie}(A,M)\stackrel{\delta}\rightarrow $$
$$\bar{H}^0_{rsym}(A,M)\rightarrow H^2_{rsym}(A,M)\rightarrow H^2_{lie}(A,M)\stackrel{\delta}\rightarrow
\bar{H}^1_{rsym}(A,M)\rightarrow\cdots$$
$$\stackrel{\delta}\rightarrow \bar{H}^{k-2}_{rsym}(A,M)\rightarrow H^k_{rsym}(A,M)
\rightarrow H^k_{lie}(A,M)\stackrel{\delta}\rightarrow \bar{H}^{k-1}_{rsym}(A,M)
\rightarrow\cdots$$

\medskip

\noindent is exact. A connected homomorphism
$$\delta: \bar{H}^k_{lie}(A,M)\rightarrow H^{k-1}_{rsym}(A,M)$$

\noindent is induced by homomorphism
$$\delta: \bar{Z}^k_{lie}(A,M)\rightarrow Z^{k-1}_{rsym}(A,M),
\quad \psi\mapsto d_{rsym}\psi.$$
\end{th}

{\bf Proof.} We will check that for $k>0$ the following diagram is commutative

$$\begin{array}{ccc}
 C^k_{rsym}(A,M)&\stackrel{d_{rsym}}\longrightarrow&C^{k+1}_{rsym}(A,M)\\
\downarrow\lefteqn{f}&&\downarrow\lefteqn{f}\\
C^{k}_{lie}(A,M))&\stackrel{d_{lie}}\longrightarrow&C^{k+1}_{lie}(A,M)\\
\end{array}$$

\medskip

\noindent For $\psi\in C^k_{rsym}(A,M),$ we have

$$fd_{rsym}\psi(a_1,\ldots,a_{k+1})=$$

$$\sum_{s}(-1)^{s+k}d_{rsym}\psi(a_s,a_1,\ldots,\hat{a_s},\ldots,a_{k+1})=$$

$$\sum_{i<s}(-1)^{i+s+k+1}{\mathop{a_s\circ(\psi(a_{i},a_1,\ldots,\hat{a_i},
\ldots,\hat{a_{s}},\ldots,a_{k}))}\limline}$$
$$+\sum_{s<i}(-1)^{i+s+k}{\mathop{a_s\circ\psi(a_i,a_1,\ldots,\hat{a_s},\ldots,
\hat{a_i},\ldots,a_{k+1})}\limeq}$$

$$-\sum_{i<s}(-1)^{i+s+k+1}{\mathop{\psi(a_s\circ a_{i},a_1,\ldots,\hat{a_i},
\ldots,\hat{a_{s}},\ldots,a_{k}))}\limsim}$$
$$-\sum_{s<i}(-1)^{i+s+k}{\mathop{\psi(a_s\circ a_i,a_1,\ldots,\hat{a_s},
\ldots,\hat{a_i}.\ldots,a_{k+1})}\limsim}$$

$$+\sum_{i<j<s}(-1)^{s+i+k}{\mathop{\psi(a_s,a_1,\ldots,\hat{a_i},\ldots,a_{j-1},
[a_i,a_j],\ldots,\hat{a_s},\ldots,a_{k+1})}\limsimeq}$$

$$+\sum_{i<s<j}(-1)^{s+i+k}{\mathop{\psi(a_s,a_1,\ldots,\hat{a_i},\ldots,
\hat{a_s},\ldots,a_{j-1},[a_i,a_j],\ldots,a_{k+1})}\limapprox}$$

$$+\sum_{s<i<j}(-1)^{s+i+k+1}{\mathop{\psi(a_s,a_1,\ldots,\hat{a_s},\ldots,
\hat{a_i},\ldots,a_{j-1},[a_i,a_j],\ldots,a_{k+1})}\limequiv}$$

$$-\sum_{i<s}(-1)^{i+s+k+1}{\mathop{\psi(a_s,a_1,\ldots,\hat{a_i},
\ldots,\hat{a_{s}},\ldots,a_{k}))\circ a_i}\limline}$$
$$-\sum_{s<i}(-1)^{i+s+k}{\mathop{\psi(a_s,a_1,\ldots,\hat{a_s},\ldots,
\hat{a_i}.\ldots,a_{k+1})\circ a_i}\limeq},$$

\noindent and

$$d_{lie}f\psi(a_1,\ldots,a_{k+1})=$$

$$\sum_{i<j}(-1)^if\psi(a_1,\ldots,\hat{a_i},\ldots,a_{j-1},[a_i,a_j],
\ldots,a_{k+1})$$
$$+\sum_i(-1)^i[f\psi(a_1,\ldots,\hat{a_i},
\ldots,a_{k+1}),a_i]=$$

$$\sum_{s<i<j}(-1)^{i+s+k+1}{\mathop{\psi(a_s,a_1,\ldots,\hat{a_s},\ldots,
\hat{a_i},\ldots,a_{j-1},[a_i,a_j],\ldots,a_{k+1})}\limequiv}$$

$$+\sum_{i<s<j}(-1)^{i+k+s}{\mathop{\psi(a_s,a_1,\ldots,\hat{a_i},\ldots,\hat{a_s},\ldots,
a_{j-1},[a_i,a_j],\ldots,a_{k+1})}\limapprox}$$

$$+(-1)^{i+j+k}{\mathop{\psi([a_i,a_j],a_1,\ldots,\hat{a_i},\ldots,\hat{a_j},\ldots,a_{k+1})}
\limsim}$$

$$+\sum_{i<j<s}(-1)^{i+k+s}{\mathop{\psi(a_s,a_1,\ldots,\hat{a_i},\ldots,a_{j-1},[a_i,a_j],
\ldots,\hat{a_s},\ldots,a_{k+1})}\limsimeq}$$

$$+\sum_{s<i}(-1)^{i+s+k}{\mathop{[\psi(a_s,a_1,\ldots,\hat{a_s},\ldots,
\hat{a_i},\ldots,a_{k+1}),a_i]}\limeq}$$

$$+\sum_{i<s}(-1)^{i+s+k+1}{\mathop{[\psi(a_s,a_1,\ldots,\hat{a_i},\ldots,\hat{a_{s}},
\ldots,a_{k+1}),a_i]}\limline}$$

\noindent Thus, according to right-symmetric identity,
$$fd_{rsym}\psi=d_{lie}f\psi, \quad \forall \psi\in C^k_{rsym}(A,M),\quad\forall k>0.$$

\noindent So, a short exact sequence of cochain complexes takes place
$$0\rightarrow \oplus_{k>0}\bar{C}^k_{rsym}(A,M)\rightarrow
\oplus_{k>0}C^{k}_{rsym}(A,M)\rightarrow \oplus_{k>0}C^k_{lie}(A,M)\rightarrow 0.$$

\noindent In particular, a long cohomological sequence
$$\bar{H}^0_{rsym}(A,M)\rightarrow H^2_{rsym}(A,M)\rightarrow H^2_{lie}(A,M)\rightarrow\cdots$$
$$\rightarrow \bar{H}^{k-2}_{rsym}(A,M)\rightarrow H^k_{rsym}(A,M)\rightarrow
H^k_{lie}(A,M)\rightarrow\cdots$$

\medskip

\noindent is exact. The exactness of the beginning part
$$0\rightarrow Z^1_{rsym}(A,M)\rightarrow Z^1_{lie}(A,M)\rightarrow
\bar{H}^0_{rsym}(A,M)\rightarrow H^2_{rsym}(A,M)$$

\noindent we check directly. It is clear that $d_{lie}\psi=0,$ if
$d_{rsym}\psi=0, \; \psi\in C^1(A,M)$. So, the natural homomophism
$Z^1_{rsym}(A,M)\rightarrow Z^1_{lie}(A,M)$ is a monomorphism.
Let $\delta\phi, \; \phi\in Z^1_{lie}(A,M),$  gives us
a trivial class in $\bar{H}^0(A,M)=\bar{Z}^0_{rsym}(A,M).$
Then $\phi\in Z^1_{rsym}(A,M),$ since $\delta\phi=d_{rsym}\phi.$
Suppose that $\sigma\in \bar{Z}^0_{rsym}(A,M)$ is a coboundary in
$Z^2_{rsym}(A,M),$ say $\sigma=d_{rsym}\omega,$ for some
$\omega\in C^1_{rsym}(A,M).$ Then $d_{rsym}\omega(a,b)=\sigma(a,b)=
\sigma(b,a)=d_{rsym}\omega(b,a),$ for any $a,b\in A.$ This means
that $d_{lie}\omega=0.$

The theorem is proved completely. $\bullet$

\subsection{Cup product in right-symmetric cohomologies.}\label{cup}
\begin{th}
Assume that a cup product of $A-$modules $\cup: M\times N\rightarrow S$
is given. Then a bilinear map
$$C^{*+1}_{rsym}(A,M)\times C^*_{lie}(A,N)\rightarrow C^{*+1}_{rsym}(A,S), \;
(\psi,\phi)\mapsto \psi\cup \phi,$$

\noindent defined by

$$C^{k+1}_{rsym}(A,M)\times C^l_{lie}(A,N)\rightarrow C^{k+l+1}_{rsym}(A,S), \;
(\psi, \phi)\mapsto \psi\cup \phi,$$

$$\psi\cup\phi(a_0,a_1,\ldots,a_{k+l})=$$

$$\sum_{\begin{array}{c}
\sigma\in Sym_{k+l},\\
\sigma(1)<\cdots <\sigma(k),\\
\sigma(k+1)<\cdots <\sigma(k+l) \end{array}}
sgn\,\sigma\,
\psi(a_0,a_{\sigma(1)},\ldots,a_{\sigma(k)})\cup \phi(a_{\sigma(k+1)},\ldots,
a_{\sigma(k+l)}).
$$

\noindent is also a cup product:
$$(\alpha\cup\beta){\rho_{rsym}}(x)=\alpha{\rho_{rsym}}(x)\cup\beta
+\alpha\cup \beta\rho_{lie}(x),\;\;\forall\alpha\in
C^{*+1}_{rsym}(A,M),\;\forall\beta\in C^*_{lie}(A,N).$$

\noindent Moreover,
\begin{equation}
d_{rsym}(\psi\cup\phi)=d_{rsym}\psi\cup\phi-(-1)^k\psi\cup d_{lie}\phi,
\label{cupcoh}
\end{equation}

\noindent for any $\psi\in C^{k+1}_{rsym}(A,M), \; \phi\in C^l_{lie}(A,N),
\;\; k,l\ge 0.$
\end{th}

{\bf Proof.} By theorem~\ref{chev}
$$F(\eta\cup\phi)=F(\eta)\cup\phi,$$
$$F(\alpha\rho_{lie}(x))=(F\alpha)\rho_{rsym}(x),$$

\noindent for any $\eta\in C^k_{lie}(A,C^1(A,M)),\;\alpha\in
C^{k+1}_{rsym}(A,M), \; \phi\in C^l_{lie}(A,N), \;x\in A.$

Prolongate the cup product $M\times N\rightarrow S,\; (m,n)\mapsto m\cup n,$
of $A-$ modules to a cup product of $A^{lie}$-modules
$$C^1(A,M)^{lie}\times N^{lie}\rightarrow C^1(A,S)^{lie},$$
$$(f,n)\mapsto f\cup n,\;\; (f\cup n)(a)=f(a)\cup n.$$

\noindent Check the correctness of this definition:
$$((f\cup n)\circ(a))(b)=$$
$$d_{rsym}(f\cup n)(b,a)=$$
$$b\circ ((f\cup n)(a))-(f\cup n)(b\circ a)+((f\cup n)(b))\circ a=$$
$$b\circ (f(a)\cup n)-f(b\circ a)\cup n+(f(b)\cup n)\circ a=$$
$$(b\circ f(a))\cup n-f(b\circ a)\cup n+(f(b)\circ a)\cup n+f(b)\cup [n,a]=$$
$$(d_{rsym}f(b,a))\cup n+f(b)\cup [n,a]=$$
$$(f\circ a)(b)\cup n+f(b)\cup [n,a]=$$
$$((f\circ a)\cup n+(f\cup [n,a]))(b).$$

Thus, we have a cup product of Chevalley-Eilenberg cochain complexes~\cite{Serrehoch}
$$C^k_{lie}(A,C^1(A,M))\times C^l_{lie}(A,N)\rightarrow C^{k+l}_{lie}(A,C^1(A,S))$$
$$\{(\eta\cup\phi)(a_1,\ldots,a_{k+l})\}(a_0)=$$
$$\sum_{\begin{array}{c}
\sigma\in Sym_{k+l},\\
\sigma(1)<\cdots <\sigma(k),\\
\sigma(k+1)<\cdots <\sigma(k+l) \end{array}}
sgn\,\sigma\,
\{\eta(a_{\sigma(1)},\ldots,a_{\sigma(k)})\cup \phi(a_{\sigma(k+1)},\ldots,
a_{\sigma(k+l)})\}(a_0).
$$

\noindent We see that cup products for Chevalley-Eilenberg complexes
and right-symmetric complexes are compatible. Namely,
\begin{equation}(F\eta)\cup\phi=F(\eta\cup\phi), \label{eta}\end{equation}

\noindent for any $\eta\in C^k_{lie}(A,C^1(A,M)), \; \phi\in C^l_{lie}(A,N)$
(definition of isomorphism $F:C^k_{lie}(A,C^1(A,M))\rightarrow C^{k+1}_{rsym}(A,M)$
see~(\ref{F})). Since,~\cite{Serrehoch}
$$d_{lie}(\eta\cup\phi)=d_{lie}\eta\cup\phi+(-1)^k\eta\cup d_{lie}\phi,$$

\noindent accordingly~(\ref{eta}),
$$d_{rsym}((F\eta)\cup\phi)=d_{rsym}F(\eta\cup\phi)=$$
$$Fd_{lie}(\eta\cup\phi)=F(d_{lie}\eta\cup \phi+(-1)^k\eta\cup d_{lie}\phi)=$$
$$Fd_{lie}\eta\cup\phi+(-1)^kF\eta\cup d_{lie}\phi=$$
$$d_{rsym}F\eta\cup\phi+(-1)^kF\eta\cup d_{lie}\phi.$$

\noindent By theorem~\ref{chev} for any $\psi\in C^{k+1}_{rsym}(A,M), k\ge 0,$
there exists $\eta\in C^k(A,C^1(A,M)),$ such that $\psi=F\eta.$ Hence,
(\ref{cupcoh}) is true. $\bullet$

\begin{crl}
The cup product
$$C^{*+1}_{rsym}(A,M)\times C^*_{lie}(A,N)\rightarrow C^{*+1}_{rsym}(A,S), \;\; (\psi,\phi)\mapsto \psi\cup \phi,$$

\noindent induces a cup product of cohomology spaces
$$H^{k+1}_{rsym}(A,M)\times H^l_{lie}(A,N)\rightarrow H^{k+l+1}_{rsym}(A,S),\;\;
([\psi],[\phi])\mapsto [\psi\cup\phi],\;\; k>0, l\ge 0.$$
$$Z^1_{rsym}(A,M)\times H^l_{lie}(A,N)\rightarrow H^{l+1}_{rsym}(A,S),
\;\; (\psi, [\phi])\mapsto [\psi\cup\phi],\; l\ge 0.$$
\end{crl}

{\bf Proof.}
$$Z^{k+1}_{rsym}(A,M)\cup Z^l_{lie}(A,N)\subseteq Z^{k+l+1}_{rsym}(A,S),
\;\; k,l\ge 0,$$

$$B^{k+1}_{rsym}(A,M)\cup Z^l_{lie}(A,N)\subseteq B^{k+l+1}_{rsym}(A,S),
\;\; k>0, \;l\ge 0,$$

$$Z^{k+1}_{rsym}(A,N)\cup B^{l}_{lie}(A,N)\subseteq B^{k+l+1}_{rsym}(A,S),
\;\; k,l\ge 0.\bullet$$

Notice that for any module $M$ of right-symmetric algebra $A$ and trivial
$A$-module ${\cal K}$ there exists a natural cup product
$$M\times {\cal K}\rightarrow M,\;\; (m,\lambda)\mapsto m\lambda.$$

\noindent So, we have a pairing of cohomology spaces
$$H^*_{rsym}(A,M)\times H^*_{lie}(A,{\cal K})\rightarrow H^*_{rsym}(A,M).$$

\noindent In particular, $H^*_{rsym}(A,M)$ has a natural structure of module over
$H^*_{lie}(A,{\cal K}).$ As it turned out in some cases $H^*_{rsym}(A,M)$ is a free
$H^*_{lie}(A,{\cal K})$-module. In section~\ref{calc} we will see that this 
is the case, if $A=gl_n^{rsym}.$

Denote by $\bar M$ antisymmetric $A$-module obtained from $M$ by
${\bar r}_a=r_a-l_a, {\bar l}_a=0.$ One can construct another cup product
$${\cal K}\times M\rightarrow {\bar M},\;\; \lambda\cup m=\lambda m.$$

\noindent We use this cup product in consideration of right-symmetric cohomologies
for $A=W_n^{rsym},$ section~\ref{calc}.

\section{Deformations of right-symmetric algebras.}

\subsection{Deformation equations}
We will follow the Gerstenhaber theory of  deformations of algebras \cite{Gerst}. 
Let $A$ be a right-symmetry algebra over a field ${\cal K}$ of any characteristic $p.$
Let ${\cal K}((t))$ be a fraction field for formal power series algebra ${\cal K}[[x]].$
Extend the main field ${\cal K}$ until ${\cal K}((t))$ and construct on the vector space
$A\otimes {\cal K}((t))$ a new right-symmetric multiplication
$$\mu_t=\mu_0+t\mu_1+t^2\mu_2+\ldots,$$

\noindent where
$$\mu_i\in C^2_{rsym}(A,M), \quad i=0,1,2,\ldots, \mbox{\;and \;} \mu_0(a,b)=a\circ b.$$

The right-symmetric condition for $\mu_t$ in terms of $\mu_k$ can be regarded
as the following {\it deformation equations}
$$\mu_1\in Z^2_{rsym}(A,A), \eqno{(DFR.1)} $$
$$\sum_{l=1}^{k-1}\mu_l\star\mu_{k-l}=-d_{rsym}\mu_k,  \eqno{(DFR.k)}$$
$$k=2,3,\ldots ,$$

\noindent where
$$(\psi\star\phi)(a,b,c)=
\psi(a,\phi(b,c))-\psi(\phi(a,b),c)-\psi(a,\phi(c,b))+\psi(\phi(a,c),b),$$

$\psi, \phi\in C^2_{rsym}(A,M).$

Right-symmetric deformations $\mu_t,$ and $\nu_t$ are said to be {\it equivalent,}
if there exists a map
$$g_t=g_0+tg_1+tg_2+\cdots, \quad g_k\in C^1_{rsym}(A,A), \;\; k=0,1,2,\ldots ,$$

\noindent with an identity map $g_0,$ such that
$$g_t^{-1}(\mu_t(g_t(a),g_t(b)))=\nu_t(a,b),  \; \forall a,b\in A.$$

\noindent In particular, for equivalent deformations $\mu_t, \nu_t,$ should be
$$\nu_1=\mu_1+d_{rsym}g_1.$$

\noindent In other words the first deformation terms, so called {\it local deformations}
will define equivalent 2-right-symmetry cohomology classes $[\mu_1]=[\nu_1].$

Converse, suppose that there is given a 2-cocycle of right-symmetric algebra with
coefficients in the regular module, $\psi\in Z^2_{rsym}(A,A),$ with a cohomology class
$[\psi]\in H^2_{rsym}(A,A).$ One can take $\mu_1:=\psi,$ and try to construct
$\mu_k$ that will satisfy deformation equations. Evidently, (DFR.1) is true.
We will say that {\it local deformation $\mu_1=\psi$ can be prolongated to a global
deformation until $k$-th term,} if there exist $\mu_2,\ldots,\mu_k,$ such that
equations (DFR.k) are true. If this is the case for any $k>0,$ we will say
that {\it local deformation $\mu_1$ can be prolongated until global deformation $\mu_t$}
or, equivalently, that $\mu_t$ is global deformation or prolongation of $\mu_1.$ Set,
$$Obs_k(\psi)=\sum_{l=1}^{k-1}\mu_l\star\mu_{k-l}.$$

\noindent Notice that the definition of $Obs_k(\psi)$ depends not only from $\psi$ but, also
from the first $k-1$ terms of deformation.

\subsection{Third cohomologies as obstruction}
\begin{prp} Suppose that a local deformation $\mu_1=\psi$ can be prolongated
to a global deformation until $(k-1)$-th term. Then, $Obs_k(\psi)\in Z^3_{rsym}(A,A)$
and the prolongation of $\mu_1$ until $k$-th term is possible, if and only if
$[Obs_k(\psi)]=0.$
\end{prp}

{\bf Proof.} For $\alpha\in C^{k+1}(A,A), \beta\in C^{l+1}(A,A)$ define
multiplications $\alpha\ast\beta\in C^{k+l+1}(A,A), \;\; \alpha\smile\beta\in
C^{k+l+2}(A,A)$ by
$$\alpha\ast\beta(a_1,\ldots,a_{k+l+1})=$$
$$\sum_{s=1}^{k+1}(-1)^{(s+1)l}\alpha(a_1,\ldots,a_{s-1},\beta(a_{s+1},\ldots,a_{s+l}),
a_{s+l+1},\ldots,a_{k+l+1}).$$

$$\alpha\smile\beta(a_1,\ldots,a_{k+l+2})=$$
$$\alpha(a_1,\ldots,a_{k+1})\circ\beta(a_{k+2},\ldots,a_{k+l+2}).$$

\noindent Then,
$$\psi\star\phi(a,b,c)=\psi\ast\phi(a,c,b)-\psi\ast\phi(a,b,c), \;\; \psi,\phi\in C^2(A,A),$$

\noindent and
$$d_{rsym}Obs_k(\psi)(a_0,a_1,a_2,a_3)=$$
$$\sum_{l+s=k, l>0, s>0}\sum_{\sigma\in Sym_3}
sgn\,\sigma\,d_{ass}(\mu_l\ast\mu_{s})(a_0,a_{\sigma(1)},a_{\sigma(2)}, a_{\sigma(3)})$$

\noindent where $d_{ass}$ means Hochshild coboundary operator as in associative algberas.
By [\cite{Gerst1}, \S 7, Th.3],
$$d_{ass}\alpha\ast\beta=\alpha\ast d_{ass}\beta-d_{ass}\alpha\ast\beta
-\alpha\smile\beta+\beta\smile\alpha, \;\; \alpha, \beta\in C^2(A,A).$$

\noindent Notice that
$$\sum_{l+s=k, l>0, s>0}{\mu_l\smile\mu_s-\mu_s\smile\mu_l}=
\sum_{l+s-k, l>0, s>0}\mu_l\smile\mu_s-\sum_{l+s-k, l>0, s>0}\mu_l\smile\mu_s=0.$$

\noindent Hence, according to conditions (DFR.$l$), $l<k,$

$$d_{rsym}Obs_k(\psi)(a_0,a_1,a_2,a_3)=$$

$$\sum_{l+s=k, l>0, s>0} \sum_{\sigma\in Sym_3}
sgn\,\sigma \,d_{ass}(\mu_l\smile \mu_s)
(a_0,a_{\sigma(1)},a_{\sigma(2)},a_{\sigma(3)})=$$

$$\sum_{l+s=k, l>0, s>0} \sum_{\sigma\in Sym_3}
sgn\,\sigma\,\mu_l\ast d_{ass}\mu_s(a_0,a_{\sigma(1)},a_{\sigma(2)},a_{\sigma(3)})$$
$$\;\;-sgn\,\sigma\,d_{ass}\mu_l\ast\mu_s(a_0,a_{\sigma(1)},a_{\sigma(2)},a_{\sigma(3)})=$$

$$\sum_{l+s=k, l>0, s>0} \sum_{\sigma\in Sym_3}
sgn\,\sigma\,\mu_l(d_{ass}\mu_s(a_0,a_{\sigma(1)},a_{\sigma(2)}),a_{\sigma(3)})$$
$$-sgn\,\sigma\,\mu_l(a_0, d_{ass}\mu_s(a_{\sigma(1)},a_{\sigma(2)},a_{\sigma(3)}))$$
$$\;\;-sgn\,\sigma\,d_{ass}\mu_l(\mu_s(a_0,a_{\sigma(1)}),a_{\sigma(2)},a_{\sigma(3)})$$
$$+sgn\,\sigma\,d_{ass}\mu_l(a_0,\mu_s(a_{\sigma(1)}, a_{\sigma(2)}),a_{\sigma(3)})$$
$$-sgn\,\sigma\,d_{ass}\mu_l(a_0,a_{\sigma(1)},\mu_s(a_{\sigma(2)},a_{\sigma(3)}))=$$

$$\sum_{l+s=k, l>0, s>0}$$
$$\{\mu_l(d_{rsym}\mu_s(a_0,a_1,a_2),a_3)
-\mu_l(d_{rsym}\mu_s(a_0,a_1,a_3),a_2)
+\mu_l(d_{rsym}\mu_s(a_0,a_2,a_3),a_1)$$
$$-\mu_l(a_0, d_{rsym}\mu_s(a_{1},a_{2},a_{3}))
+\mu_l(a_0, d_{rsym}\mu_s(a_{2},a_{1},a_{3}))
-\mu_l(a_0, d_{rsym}\mu_s(a_{3},a_{1},a_{2}))$$
$$-d_{rsym}\mu_l(\mu_s(a_0,a_{1}),a_{2},a_{3})
+d_{rsym}\mu_l(\mu_s(a_0,a_{2}),a_{1},a_{3})
-d_{rsym}\mu_l(\mu_s(a_0,a_{3}),a_{1},a_{2})$$
$$+d_{rsym}\mu_l(a_0,\mu_s(a_{1}, a_{2}),a_{3})
-d_{rsym}\mu_l(a_0,\mu_s(a_{2}, a_{1}),a_{3})$$
$$-d_{rsym}\mu_l(a_0,\mu_s(a_{1}, a_{3}),a_{2})
+d_{rsym}\mu_l(a_0,\mu_s(a_{3}, a_{1}),a_{2})$$
$$+d_{rsym}\mu_l(a_0,\mu_s(a_{2}, a_{3}),a_{1})
-d_{rsym}\mu_l(a_0,\mu_s(a_{3}, a_{2}),a_{1})\}=$$
$$S_1+S_2,$$

\noindent where

$$S_1=\sum_{l+s=k, l>0, s>0}\;\sum_{s_1+s_2=s, s_1, s_2>0}$$
$$\{-\mu_l(\mu_{s_1}\star\mu_{s_2}(a_0,a_1,a_2),a_3)
+\mu_l(\mu_{s_1}\star\mu_{s_2}(a_0,a_1,a_3),a_2)
-\mu_l(\mu_{s_1}\star\mu_{s_2}(a_0,a_2,a_3),a_1)$$
$$+\mu_l(a_0, \mu_{s_1}\star\mu_{s_2}(a_{1},a_{2},a_{3}))
-\mu_l(a_0,\mu_{s_1}\star\mu_{s_2}(a_{2},a_{1},a_{3}))
+\mu_l(a_0, \mu_{s_1}\star\mu_{s_2}(a_{3},a_{1},a_{2}))\},$$

$$S_2=\sum_{l+s=k, l>0, s>0}\;\sum_{l_1+l_2=l, l_1, l_2>0}$$
$$\{\mu_{l_1}\star\mu_{l_2}(\mu_s(a_0,a_{1}),a_{2},a_{3})
-\mu_{l_1}\star\mu_{l_2}(\mu_s(a_0,a_{2}),a_{1},a_{3})
+\mu_{l_1}\star\mu_{l_2}(\mu_s(a_0,a_{3}),a_{1},a_{2})$$
$$-\mu_{l_1}\star\mu_{l_2}(a_0,\mu_s(a_{1}, a_{2}),a_{3})
+\mu_{l_1}\star\mu_{l_2}(a_0,\mu_s(a_{2}, a_{1}),a_{3})$$
$$+\mu_{l_1}\star\mu_{l_2}(a_0,\mu_s(a_{1}, a_{3}),a_{2})
-\mu_{l_1}\star\mu_{l_2}(a_0,\mu_s(a_{3}, a_{1}),a_{2})$$
$$-\mu_{l_1}\star\mu_{l_2}(a_0,\mu_s(a_{2}, a_{3}),a_{1})
+\mu_{l_1}\star\mu_{l_2}(a_0,\mu_s(a_{3}, a_{2}),a_{1})\}.$$

We have

$$S_1=\sum_{l+s=k, l>0, s>0}\;\sum_{s_1+s_2=s, s_1>0, s_2>0}$$

$$\{-\mathop{\mu_l(\mu_{s_1}\ast\mu_{s_2}(a_0,a_1,a_2),a_3)}\limline
+\mathop{\mu_l(\mu_{s_1}\ast\mu_{s_2}(a_0,a_1,a_3),a_2)}\limeq
-\mathop{\mu_l(\mu_{s_1}\ast\mu_{s_2}(a_0,a_2,a_3),a_1)}\limli$$
$$+\mathop{\mu_l(a_0, \mu_{s_1}\ast\mu_{s_2}(a_{1},a_{2},a_{3}))}\limline
-\mathop{\mu_l(a_0,\mu_{s_1}\ast\mu_{s_2}(a_{2},a_{1},a_{3}))}\limsim
+\mathop{\mu_l(a_0, \mu_{s_1}\ast\mu_{s_2}(a_{3},a_{1},a_{2}))}\limsimeq$$

$$+\mathop{\mu_l(\mu_{s_1}\ast\mu_{s_2}(a_0,a_2,a_1),a_3)}\limsim
-\mathop{\mu_l(\mu_{s_1}\ast\mu_{s_2}(a_0,a_3,a_1),a_2)}\limsimeq
+\mathop{\mu_l(\mu_{s_1}\ast\mu_{s_2}(a_0,a_3,a_2),a_1)}\limsmile$$
$$-\mathop{\mu_l(a_0, \mu_{s_1}\ast\mu_{s_2}(a_{1},a_{3},a_{2}))}\limeq
+\mathop{\mu_l(a_0,\mu_{s_1}\ast\mu_{s_2}(a_{2},a_{3},a_{1}))}\limli
-\mathop{\mu_l(a_0, \mu_{s_1}\ast\mu_{s_2}(a_{3},a_{2},a_{1}))}\limsmile\}=$$

$$\sum_{l+s_1+s_2=k, l>0, s_1>0, s_2>0}\;\;\sum_{\sigma\in Sym_3}
-sgn\,\sigma\,\mu_{l}\ast(\mu_{s_1}\ast\mu_{s_2})(a_0,a_{\sigma(1)},
a_{\sigma(2)},a_{\sigma(3)}),$$

\noindent and

$$S_2=\sum_{l+s=k, l>0, s>0}\;\,\sum_{l_1+l_2=l, l_1>0, l_2>0}$$

$$\{\mathop{\mu_{l_1}\ast\mu_{l_2}(\mu_s(a_0,a_{1}),a_{2},a_{3})}\limline
-\mathop{\mu_{l_1}\ast\mu_{l_2}(\mu_s(a_0,a_{2}),a_{1},a_{3})}\limli
+\mathop{\mu_{l_1}\ast\mu_{l_2}(\mu_s(a_0,a_{3}),a_{1},a_{2})}\limsim$$
$$-\mathop{\mu_{l_1}\ast\mu_{l_2}(a_0,\mu_s(a_{1}, a_{2}),a_{3})}\limline
+\mathop{\mu_{l_1}\ast\mu_{l_2}(a_0,\mu_s(a_{2}, a_{1}),a_{3})}\limli$$
$$+\mathop{\mu_{l_1}\ast\mu_{l_2}(a_0,\mu_s(a_{1}, a_{3}),a_{2})}\limeq
-\mathop{\mu_{l_1}\ast\mu_{l_2}(a_0,\mu_s(a_{3}, a_{1}),a_{2})}\limsim$$
$$-\mathop{\mu_{l_1}\ast\mu_{l_2}(a_0,\mu_s(a_{2}, a_{3}),a_{1})}\limsimeq
+\mathop{\mu_{l_1}\ast\mu_{l_2}(a_0,\mu_s(a_{3}, a_{2}),a_{1})}\limcong$$

$$-\mathop{\mu_{l_1}\ast\mu_{l_2}(\mu_s(a_0,a_{1}),a_{3},a_{2})}\limeq
+\mathop{\mu_{l_1}\ast\mu_{l_2}(\mu_s(a_0,a_{2}),a_{3},a_{1})}\limsimeq
-\mathop{\mu_{l_1}\ast\mu_{l_2}(\mu_s(a_0,a_{3}),a_{2},a_{1})}\limcong$$
$$+\mathop{\mu_{l_1}\ast\mu_{l_2}(a_0,a_3, \mu_s(a_{1}, a_{2}))}\limsim
-\mathop{\mu_{l_1}\ast\mu_{l_2}(a_0,a_3,\mu_s(a_{2}, a_{1}))}\limcong$$
$$-\mathop{\mu_{l_1}\ast\mu_{l_2}(a_0,a_2,\mu_s(a_{1}, a_{3}))}\limli
+\mathop{\mu_{l_1}\ast\mu_{l_2}(a_0,a_2,\mu_s(a_{3}, a_{1}))}\limsimeq$$
$$+\mathop{\mu_{l_1}\ast\mu_{l_2}(a_0,a_1,\mu_s(a_{2}, a_{3}))}\limline
-\mathop{\mu_{l_1}\ast\mu_{l_2}(a_0, a_1,\mu_s(a_{3}, a_{2}))}\limeq\}=$$

$$\sum_{l_1+l_2+s=k, l_1>0, l_2>0, s>0}sgn\,\sigma\,(\mu_{l_1}\ast\mu_{l_2})\ast\mu_s
(a_0,a_{\sigma(1)},a_{\sigma(2)},a_{\sigma(3)}),$$

Let $\alpha,\beta,\gamma\in C^2(A,A).$ Then,

$$\{\alpha\ast(\beta\ast\gamma)-(\alpha\ast\beta)\ast\gamma\}(a,b,c,d)=$$

$$\alpha(\beta\ast\gamma(a,b,c),d)+\alpha(a,\beta\ast\gamma(b,c,d))$$
$$-\alpha\ast\beta(\gamma(a,b),c,d)+\alpha\ast\beta(a,\gamma(b,c),d)
-\alpha\ast\beta(a,b,\gamma(c,d))=$$

$$\mathop{\alpha(\beta(\gamma(a,b),c),d)}\limline
-\mathop{\alpha(\beta(a,\gamma(b,c)),d)}\limcong$$
$$+\mathop{\alpha(a,\beta(\gamma(b,c),d))}\limeq
-\mathop{\alpha(a,\beta(b,\gamma(c,d)))}\limsim$$

$$-\mathop{\alpha(\beta(\gamma(a,b),c),d)}\limline
+\alpha(\gamma(a,b),\beta(c,d))$$
$$+\mathop{\alpha(\beta(a,\gamma(b,c)),d)}\limcong
-\mathop{\alpha(a,\beta(\gamma(b,c),d))}\limeq$$
$$-\alpha(\beta(a,b),\gamma(c,d))+
\mathop{\alpha(a,\beta(b,\gamma(c,d)))}\limsim= $$

$$-\alpha(\beta(a,b),\gamma(c,d))+\alpha(\gamma(a,b),\beta(c,d))=$$

So, for any $\alpha,\beta,\gamma\in C^2(A,A),$
$$\alpha\ast(\beta\ast\gamma+\gamma\ast\beta)-(\alpha\ast\beta)\ast\gamma-(\alpha\ast\gamma)\ast\beta=0$$

\noindent For these reasons,
$$S_1=$$
$$\sum_{s_1+s_2+s_3=k, s_1>0, s_2>0, s_3>0}\;\;\sum_{\sigma\in Sym_3}
-sgn\,\sigma\,\mu_{s_1}\ast(\mu_{s_2}\ast\mu_{s_3})(a_0,a_{\sigma(1)},
a_{\sigma(2)},a_{\sigma(3)})=$$
$$-\sum_{l_1+l_2+l_3=k, l_1>0, l_2>0, l_3>0}sgn\,
\sigma\,(\mu_{l_1}\ast\mu_{l_2})\ast\mu_{l_3}
(a_0,a_{\sigma(1)},a_{\sigma(2)},a_{\sigma(3)})=$$
$$=-S_2.$$

So, we prove that $d_{rsym}Obs_k(\psi)=0,$ if $d_{rsym}Obs_l(\psi)=0,$
for any $0<l<k. \quad \bullet$

\begin{crl}\label{prod} If $H^3_{rsym}(A,A)=0,$ then any local deformation
can be prolongated.
\end{crl}

\subsection{Steenrod squares\label{Sq}}

Let $char\,k=p>0.$ In this subsection we recall Gerstenhaber's construction
\cite{Gerst} of the homomorpism
$$Sq: Z^1_{rsym}(A,A)\rightarrow Z^2_{rsym}(A,A), \quad D\mapsto Sq D$$

\noindent regarding the right-symmetric algebras.
For any derivation $D\in Z^1_{rsym}(A,A)$ its $p$-th power $D^p$
is also a derivation, $D^p\in Z^1_{rsym}(A,A).$
The proof is based on the following property of binomial coefficients:
an integer ${p\choose a}$ can be divided into $p,$ if $0<a<p.$ Then,
$$D^p(a\circ b)-D^p(a)\circ b-a\circ D^p(b)=$$
$$\sum_{i=1}^{p-1}{p\choose i}D^i(a)\circ D^{p-i}(b)\equiv 0(mod\,p).$$

\noindent In particular, we can consider integers  ${p\choose i}/p,\;\; 0<i<p$
by modulus $p$ and introduce 2-cocycle $Sq D$ with coefficients in the regular
module
$$Sq\, D(a,b)=\sum_{l=1}^{p-1}D^i(a)\circ D^{p-i}(b)/{i!(p-i)!}.$$

\noindent This cocycle is called the Steenrod Square of derivation $D$
and can be interpreted as an obstruction to prolongation of derivation
to automorphism.

\section{Calculations \label{calc}}

\subsection{Standard 2-cocycles of right-symmetric algebras.\label{standard}}

In this subsection we will give the second interpretation of the identity
$d^3_{rsym}m=0,\; m\in M,$ mentioned in section~\ref{cohom}.

\begin{prp} i) Let $\tilde M$ be a module over right-symmetric algebra $A$
and $M$ is its submodule. Suppose that for $\tilde m\in \tilde M,$
$$(\tilde m,a,b)\in M,\;\; \forall a,b\in A.$$

\noindent Then 2-cochain $\psi_{\tilde m}\in C^2(A,M)$ defined by
$$\psi_{\tilde m}(a,b)={\tilde m}\circ(a\circ b)-({\tilde m}\circ a)\circ b,$$

\noindent is symmetric 2-cocycle, $\psi_{\tilde m}\in {\bar Z}^2(A,M).$

If $\tilde m\circ a\in M,\;\forall a\in A,$ then
$[\psi_{\tilde m}]=[\phi_{\tilde m}],$ where
$$\phi_{\tilde m}(a,b)=a\circ(\tilde m\circ b).$$
\end{prp}

Notice that in the denotions of section~\ref{nabla}, $\psi_{\tilde m}=\nabla(\tilde m).$

{\bf Proof.} $$\psi_{\tilde m}(a,b)=(m,a,b)=(m,b,a)=\psi_{\tilde m}(b,a),$$

$$d_{rsym}\psi_{\tilde m}(a,b,c)=$$
$$a\circ(\tilde m,b,c)-a\circ(\tilde m,c,b)-({\tilde m},a\circ b,c)+({\tilde m},a\circ c,b)+({\tilde m},a,[b,c])-({\tilde m},a,b)\circ c+({\tilde m},a,c)\circ b=$$
$$-({\tilde m},a\circ b,c)+({\tilde m},a,b\circ c)-({\tilde m},a,b)\circ c$$
$$+({\tilde m},a\circ c,b)-({\tilde m},a,c\circ b)+({\tilde m},a,c)\circ b=$$
$${\tilde m}\circ (a,b,c)-({\tilde m}\circ a,b,c)$$
$$-{\tilde m}\circ (a,c,b)+({\tilde m}\circ a,c,b)=$$
$$0.$$

If $\tilde m\circ a\in M,$ then we can introduce a linear map
$\omega:A\rightarrow M, \; a\mapsto \tilde m\circ a.$ We obtain
$$\psi_{\tilde m}(a,b)+d_{rsym}\omega(a,b)=$$
$$\psi_{\tilde m}(a,b)-\omega(a\circ b)+\omega(a)\circ b+a\circ\omega(b)=$$
$$\psi_{\tilde m}(a,b)-\tilde m\circ(a\circ b)+(\tilde m \circ a)\circ b+
a\circ(\tilde m\circ b)=$$
$$a\circ(\tilde m\circ b).$$

\noindent In other words, $\phi_{\tilde m}=\psi_{\tilde m}+d_{rsym}
\omega.\bullet$

\subsection{Semi-center and derivations of $W_n^{rsym}$ \label{semi-center}}
Suppose that  $A$ is an algebra with multiplications $(a,b)\mapsto a\circ b,$
and $(a,b)\mapsto a\ast b,$ such that the following conditions hold
$$a\circ (b\circ c)-(a\circ b)\circ c-a\circ(c\circ b)+(a\circ c)\circ b=0,$$
$$a\ast (b\ast c)-b\ast(a\ast c)=0,$$
$$a\circ (b\ast c)-b\ast(a\circ c)=0,$$
$$(a\ast b-b\ast a- a\circ b+b\circ a)\ast c=0,$$
$$(a\circ b- b\circ a)\ast c+a\ast(c\circ b)-(a\ast c)\circ b
-b\ast(c\circ a)+(b\ast c)\circ a=0.$$

\noindent In particular, $A$ is a right-symmetric algebra.

Let $Z_l(A)$ be the left center of $A,$ $Q_l(A)$ is the left units space and
$N_{l}(A)=Z_l(A)\oplus Q_l(A)$ is the left semi-center.

\begin{th} For $A=W_n^{rsym},$ if $p=0,$ or $A=W_n({\bf m}),$ if $p>2,$
$$Z_l(A)=\{\der_i: i=1,\ldots, n\}\oplus \delta(p>0)\{\der_i^{p^{k_i}-1} :
0<k_i<m_i, i=1,\ldots,n\}\cong$$
$ {\cal K}^n\oplus \delta(p>0) {\cal
K}^{m-n},$
$$Q_l(A)=\{e=(1/n)\sum_{i=1}^nx_i\der_i\}\cong {\cal K},$$
$$Z^1_{rsym}(A,A)=\{ad\,\der_i, ad\,x_i\der_j ,
\delta(p>0)\der_i^{p^{k_i}}
: i,j,=1,\ldots,n, 0<k_i<m_i\}\cong$$
${\cal K}^n\oplus gl_n\oplus \delta(p>0) {\cal K}^{m-n},$

\medskip

\noindent where $m=\sum_{i=1}^nm_i,$ and $\delta(p>0)=0,$ if $p=0$
and =1, if $p>0.$

\end{th}

{\bf Proof.} Any derivation of right-symmetric algebra $A$
induces a derivation of Lie algebra $A^{lie}:$
$$d_{rsym}(a,b)=0,\forall a,b\in A\Rightarrow d_{lie}f(a,b)=
d_{rsym}f(a,b)-d_{rsym}(b,a)=0, \forall a,b\in A.$$

\noindent The corresponding homomorphism
$Z^1_{rsym}(A,A)\rightarrow Z^1_{lie}(A,A)$
is monomorphism. It is known, that all Lie derivations of $W_n$ are inner,
i.e. have a form  $ad\,u\der_i,$ where $u\in U, i=1,\ldots, n.$    In case of
$W_n({\bf m}), p>0,$ outer derivations $\der_i^{p^{k_i}}, 0<k_i<m_i,
i=1,\ldots, m_i,$  appear.
Since
$$ad\, a(b\circ c)-ad\,a(b)\circ c-b\circ ad\,a(c)=$$
$$a\circ (b\circ c)-(b\circ c)\circ a-(a\circ b)\circ c+(b\circ a)\circ c-
b\circ (a\circ c)+b\circ (c\circ a)=$$
$$a\circ (b\circ c)-(a\circ b)\circ c,$$

\noindent Lie derivation $ad\,a:b\mapsto [a,b]:=a\circ b-b\circ a,$
is a right-symmetric derivation, if and only if
$$a\in A^{l.ass}.$$

\noindent It is easy to see that,
$$u\der_i\circ(v\der_j\circ w\der_s)-
(u\der_i\circ v\der_j)\circ w\der_s=\der_j\der_s(u)vw\der_i.$$

\noindent Therefore,
$$A^{l.ass}=\{u\der_i: \der_j\der_s(u)=0,\forall i,j=1,\ldots,n\}=
\{\der_i, x_i\der_j: i,j=1,\ldots,n\}.$$

In case of $p>0,$ direct calculations show that $\der_i^{p^{k_i}}\in
Z^1_{rsym}(A,A).$ Other statements of the theorem are evident. $\bullet$

\subsection{Pairing of  $W_n^{rsym}-$modules and cocycle constructions.}
\begin{th} Let $A=W_n, p=0,$ or $A=W_n({\bf m}), p>3.$
The space $H^2_{rsym}(A,A), p=0,$ has a basis
consisting of cocycle classes of four types 
$\psi^1_{s,l,r}, \psi^2_{l,r},\psi^3_{s,l}, \psi^4_{l},\;\; 
s,l,r=1,\ldots,n,$
such that
$$\psi^1_{s,l,r}(u\der_i,v\der_j) =
\delta_{j,r}x_r^{-1}(\delta_{i,s}uv\der_l -x_s\der_l(u)v\der_i),$$
$$\psi^2_{l,r}(u\der_i,v\der_j)=\delta_{j,r} x_r^{-1}\der_l(u)v\der_i,$$
$$\psi^3_{s,l}(u\der_i,v\der_j)=(\delta_{i,s}u\der_j(v)\der_l-
x_s\der_l(u)\der_j(v)\der_i),$$
$$\psi^4_l(u\der_i,v\der_j)=\der_l(u)\der_j(v)\der_i.$$

In the case of $p>3,$ the space $H^2_{rsym}(A,A)$ has a basis
with cohomological classes of
the following cocycles of five types. 
$$\psi^1_{s,l,r}(u\der_i,v\der_j) =
\delta_{j,r}x_r^{p^{m_r}-1}(\delta_{i,s}uv\der_l -x_s\der_l(u)v\der_i),$$
$$\psi^2_{l,k_l,r}(u\der_i,v\der_j)=\delta_{j,r} x_r^{p^{m_r}-1}
\der_l^{p^{k_l}}(u)v\der_i,$$
$$\psi^3_{s,l}(u\der_i,v\der_j)=(\delta_{i,s}u\der_j(v)\der_l-
x_s\der_l(u)\der_j(v)\der_i),$$
$$\psi^4_{l,k_l}(u\der_i,v\der_j)=\der_l^{p^{k_l}}(u)\der_j(v)\der_i,$$
$$\psi^5_{l,k_l}=Sq\,\der_l^{p^{k_l}},$$

\noindent where $s,l,r=1,\ldots,n,\; 0\le k_l<m_l.$
\end{th}

The proof is based on the following observations.

Suppose that $A-$module $M$ preserve the action of $N_l(A):$
$$z\circ m=0, \forall z\in Z_l(A),\quad e\circ m=m, \forall e\in Q_l(A),$$
for any $m\in M.$ Define an operator 
$$C^{k+1}_{rsym}(A,M)\rightarrow C^k_{rsym}(A,M),$$
$$i_0(a)\psi(a_1,\ldots,a_k)=\psi(a_0,a_1,\ldots,a_k)$$

\noindent and an operator

$$T: C^k_{rsym}(A,A)\rightarrow C^{k+1}_{rsym}(A,A),$$
$$T\psi(a_0,a_1,\ldots,a_k)=\sum_{i=1}^k(-1)^{i+k}a_i\ast\psi(a_0,a_1,\ldots,
\hat{a_i},\ldots,a_k).$$

\noindent Then
$$Td_{rsym}=d_{rsym}T,$$

\noindent and for any $a\in N_l(A), \psi\in Z^{k+1}_{rsym}(A,A),$
$$i_0(a)\psi\in Z^k_{lie}(A,A_{anti}).$$

Define  a pairing of regular $A-$module $A$ and antisymmetric $A-$module
$U:$
$$A\times U\rightarrow A,\quad u\der_i\cup v= uv\der_i.$$
Notice that,
$$A_{anti}\cong U\otimes Z_l(A).$$

In particular, we have pairing
$$A\times A_{anti}\rightarrow A, \;\;  u\der_i\cup v\der_j= uv\der_i.$$
Therefore we have imbedding
$$Z^1_{rsym}(A,A)\times H^k_{lie}(A,U)\rightarrow H^{k+1}_{rsym}(A,A).$$

Notice that the four types of cocycles mentioned above can be obtained from
$Z^1_{rsym}(A,A)$ (see section~\ref{semi-center}) and $H^1_{lie}(A,U),$ by 
pairing $\psi\cup\phi, \; \psi\in Z^1_{rsym}(A,A),\; \phi\in H^1_{lie}(A,U).$
Recall that $H^1_{lie}(A,U)$ can be generated by the classes of cocycles
$u\der_i\mapsto ux^{-1}_r\delta_{i,r},$ and $u\der_i\mapsto \der_i(u).$ 

Another interpretation of cocycles of types 1 and 2 can be given in terms 
of standard cocycles (see section~\ref{standard}). For simplicity consider only 
the case of $p=0.$ We have
$$\psi^1_{s,l,r}=d\omega_{s,l,r}, \mbox{for \;\;} \omega_{s,l,r}(u\der_i)=
\ln x_r[x_s\der_l,u\der_i],$$
$$\psi^2_{l,r}=d\omega_{l,r}, \mbox{for \;\;} \omega_{l,r}(u\der_i)=
\ln x_r\der_l(u)\der_i,$$

\noindent and
$$\nabla(x_s\ln x_r\der_l)(u\der_i,v\der_j)=
\der_i\der_j(x_s\ln\,x_r)uv,$$
$$\nabla(\ln x_r\der_l)(u\der_i,v\der_j)=-\delta_{i,r}\delta_{j,r}uv\der_l.$$

\noindent Therefore,
$$[\psi^1_{s,l,r}]=[\nabla(x_s\ln x_r\der_l)],$$
$$[\psi^2_{l,r}]=[\nabla(\ln x_r\der_l)],$$

\noindent because of
$$\psi^1_{s,l,r}=\nabla(x_s\ln x_r\der_l)-d^{rsym}\omega^1_{s,l,r}, \mbox
{for\;\;} \omega^1_{s,l,r}\in C^1_{rsym}(W_n,W_n), 
\;\; \omega_{s,l,r}^1(u\der_i)=\delta_{i,r}x^{-1}_rx_su\der_l,$$

$$\psi^2_{l,r}=\nabla(\ln x_r\der_l)-d^{rsym}\omega^2_{l,r}, \mbox
{for\;\;} \omega^2_{l,r}\in C^1_{rsym}(W_n,W_n), 
\;\; \omega_{l,r}^2(u\der_i)=\delta_{i,r}x^{-1}_ru\der_l.$$

\begin{th} Let $A=W_1^{rsym},$ if $p=0,$ and $A=W_1(m),$ if $p>3.$
Then $H^2_{rsym}(A,A)$ has dimension 4, if $p=0,$ 
and cohomological classes of the following cocycles generate a basic
$$\psi^1(u\der,v\der)=x^{-1}uv\der,$$
$$\psi^2(u\der,v\der)=x^{-1}\der(u)v\der,$$
$$\psi^3(u\der,v\der)=(u-x\der(u))\der(v)\der,$$
$$\psi^4(u\der,v\der)=\der(u)\der(v)\der.$$
(Recall that, $\der=\der_x,$ for $n=1.$). 

For $A=W_1(m), p>3,$ the group $H^2_{rsym}(A,A)$ is $(3m+2)-$dimensional 
and the classes of the following cocycles generate its basis
$$\psi^1(u\der,v\der)=x^{p^m-1}uv\der,$$
$$\psi^2_{k}(u\der,v\der)=x^{p^m-1}\der^{p^k}(u)v\der, \;\; 0\le k<m,$$
$$\psi^3(u\der,v\der)=(u-x\der(u))\der(v)\der,$$
$$\psi^4_{k}(u\der,v\der)=\der^{p^k}(u)\der(v)\der, \;\; 0\le k<m.$$
$$Sq\,D:(a,b)\mapsto \sum_{i=1}^{p-1} D^i(a)\circ D^{p-i}(b)/(i!(p-i)!),
\; D=\der^{p^k}, 0\le k<m.$$

\noindent Cocycles of types 1 and 2 are also Novikov cocycles.
Local deformation $\psi=\sum_{i=1}^4t_i\psi^i, p=0,$
can be prolongated if and only if $t_1t_3=0, t_2t_4=0.$ 
\end{th}

\subsection{ Right-symmetric central extensions of Novikov algebras}
Let $A$ be a Novikov algebra, $R\in Der_0 \,A:=\{D\in Der\, A:
a\circ R(b)=b\circ R(a),\;\forall a,b\in A\},$ and
$\pi:A\rightarrow \cal{K},$ a linear map, such that
$$\pi(R(a))=0, \;\; \forall a\in A.$$

Define $\psi\in C^2_{rsym}(A,{\cal K}),$  by
$$\vartheta(a,b)=\pi(a\circ R(b)).$$

\begin{lm}$\vartheta\in \bar{Z}^2_{rsym}(A,{\cal K}).$
\end{lm}

{\bf Proof.} Since, $a\circ R(b)=b\circ R(a),$ then
$$\vartheta(a,b)=\vartheta(b,a).$$

\noindent We have
$$d_{rsym}\vartheta(a,b,c)=$$

$$\pi(-\mathop{(a\circ b)\circ R(c)}\limline+\mathop{(a\circ c)\circ R(b)}\limeq
+a\circ R[b,c])=$$

$$-(a\circ R(c))\circ b-\mathop{a\circ [b,R(c)]}\limsim$$
$$+(a\circ R(b))\circ c+\mathop{a[c,R(b)]}\limsim+\mathop{a\circ R[b,c])}\limsim=$$

$$\pi(-(a\circ R(c))\circ b+(a\circ R(b))\circ c)=$$

$$\pi(-R((a\circ c)\circ b)+(R(a)\circ c)\circ b+(a\circ c)\circ R(b)$$
$$+(a\circ R(b))\circ c)=$$

$$ \pi((R(a)\circ c)\circ b+b\circ R(a\circ c)+(a\circ R(b))\circ c)=$$

$$\pi((R(a)\circ c)\circ b+b\circ (R(a)\circ c)+b\circ (a\circ R(c))
+(a\circ R(b))\circ c)=$$

$$\pi((R(a)\circ c)\circ b+\mathop{R(a)\circ (b\circ c)}\limli
+\mathop{a\circ (b\circ R(c))}\limli+(a\circ R(b))\circ c)=$$

$$\pi((R(a)\circ c)\circ b+R(a\circ (b\circ c))
-\mathop{a\circ (R(b)\circ c)}\limsimeq
+\mathop{(a\circ R(b))\circ c)}\limsimeq=$$

$$\pi(\mathop{(R(a)\circ c)\circ b}\limsmile-a\circ(c\circ R(b))
+\mathop{(a\circ c)\circ R(b))}\limsmile=$$

$$\pi(-a\circ (c\circ R(b))+R((a\circ c)\circ b)-(a\circ R(c))\circ b)=$$

$$\pi(-a\circ (c\circ R(b))-(a\circ R(c))\circ b)=$$

$$\pi(-R(a\circ (c\circ b))+R(a)\circ (c\circ b)
+\mathop{a\circ(R(c)\circ b)}\limcong
-\mathop{(a\circ R(c))\circ b)}\limcong=$$

$$\pi(R(a)\circ (c\circ b)+a\circ(b\circ R(c))-(a\circ b)\circ R(c))=$$

$$\pi(c\circ(R(a)\circ b)+a\circ(c\circ R(b))-(a\circ b)\circ R(c))=$$

$$\pi(c\circ(R(a)\circ b)+c\circ(a\circ R(b))-(a\circ b)\circ R(c))=$$

$$\pi(c\circ R(a\circ b)-(a\circ b)\circ R(c))=$$

$$0.$$

\noindent So, $\vartheta\in \bar{Z}^2_{rsym}(A,{\cal K}).\quad \bullet$

In particular, algebras
$$A=W_1^{rsym}=\{e_i: e_i\circ e_j=(i+1)e_{i+j},\;\; i,j\in {\bf Z}\}, \; p=0,$$
$$A=W_1^{rsym}(m)=\{e_i: e_i\circ e_j={i+j+1\choose i}e_{i+j}, \; -1\le i,j
\le p^m-1\}, \; p>0,$$

\noindent have right-symmetric 2-cocycles with coefficients in
the trivial $A$-module ${\cal K}.$ Prove that, the cohomological class of the
cocycle $\vartheta $ in both cases is not trivial. If $\vartheta=
d_{rsym}\omega,\; \omega\in C^1(N,{\cal K}),$ then
$$\vartheta(e_i,e_j)=-(j+1)j\delta_{i+j,-1}, \mbox{\;if\:} p=0,$$
$$\vartheta(e_i,e_j)=-(-1)^i\delta_{i+j,p^m-1}, \mbox{\;if\;} p>0,$$

$$d_{rsym}\omega(e_i,e_j)=-(i+1)\omega(e_{i+j}), \mbox{\; if\;} p=0,$$
$$d_{rsym}\omega(e_i,e_j)=-{i+j+1\choose i})\omega(e_{i+j}), \mbox{\; if\;} p>0.$$

\noindent In the case of $p=0$ we have a contradiction:
$$-2=\vartheta(e_{-2},e_1)=d_{rsym}(e_{-2},e_1)=-\omega(e_{-1}),$$
$$-6=\vartheta(e_{-3},e_2)=d_{rsym}(e_{-3},e_2)=-2\omega(e_{-1}).$$

\noindent Since, $d_{rsym}\omega(e_i,e_{p^m-i-1})=0,$
we also obtain a contradiction, if $p>0.$

\begin{th} Let $A=W_1^{rsym},$ if $p=0$ and $A=W_1(m),$ if $p>0.$
Then the second right-symmetric cohomology space $H^2_{rsym}(A,{\cal K})$ has
dimension 1  and generates by a class of cocycles
$$\vartheta(e_i,e_j)=-(j+1)j\delta_{i+j,-1}, \mbox{\;if \;\;} p=0,$$
$$\vartheta(e_i,e_j)=-(-1)^i\delta_{i+j,p^m-1}, \mbox{\;if \;\;} p>0.$$

\noindent If $A$ is considered as a Novikov algebra, then any
Novikov central extension is split: $H^2_{nov}(W_1,{\cal K})=0.$
\end{th}

Recall that Novikov cohomologies are defined in~\cite{Nov}. Central extensions 
of Cartan Type Lie algebras are described in~\cite{Dzhu1}

{\bf Proof.} For $u\in U={\cal K}[[x^{\pm 1}]]$ let $\pi(u)$ be its coefficient
at $x^{-1}.$ Then $\pi(\der(u))=0, \;\; \forall u\in U.$ 
Recall that $e_i=x^{i+1}, \;i\in {\bf Z},$ and the multiplication in $A$
is given by $a\circ b=\der(a)b, \; a,b\in U.$

Prove that an isomorphism of $A^{lie}$-modules takes place
\begin{equation}C^1(A,{\cal K})\cong U_1. \label{iso}\end{equation}

A bilinear map
$$(\;,\;): U_0\times U_1\rightarrow {\cal K},\;\; (u,v)\mapsto \pi(u\cdot v),$$

\noindent is compatible with the action of $A^{lie}:$
$$((a)_0(u),v)+(u,(a)_1(v))=$$
$$\pi(-(u\circ a)\cdot v-u\cdot(v\circ a)-u\cdot(a\circ v))=$$
$$\pi(-\der(a\cdot (u\cdot v)))=0,$$

\noindent for all $a\in A, u\in U_0,v\in U_1.$ So, we have a pairing of
$A^{lie}$-modules $(\;,\;): U_0\times U_1\rightarrow {\cal K}.$
This pairing is nondegenerate. Thus, a dual $A^{lie}$-module to $U_0$ is $U_1.$
Since,
$$(f\circ a)(b)=d_{rsym}f(b,a)=-f(b\circ a), \;\; f\in C^1(N,k), a,b\in N,$$

\noindent we see that the $A^{lie}$-module $C^1(A,{\cal K})$ is isomorphic to
the dual of $U_0.$ This ends the proving of~(\ref{iso}).

By theorem~\ref{chev}
$$H^2_{rsym}(A,{\cal K})\cong H^1_{lie}(A,C^1(A,{\cal K}))\cong 
H^1_{lie}(W_1,U_1).$$

By the resutls of Gelfand and Fuchs \cite{Fuchs} the space $H^1_{lie}(W_1,U_1)$
is 1-dimensional and generates by a class of cocycle $a\mapsto \der^2(a).$
An analogous statement is also true in the case of $p>0$ \cite{Dzhuma82}.
The corresponding right-~symmetric cocycle is the cocycle $\vartheta.$

If $\psi: A\times A\rightarrow {\cal K}$ is a cocycle for central 
extension in the category of Novikov algebras, then
$$\psi(a,b\circ c)-\psi(b,a\circ c)=0, \;\; \forall a,b,c\in A.$$

\noindent The algebra $A=W_1^{nov}$ has an element $e_0$ that has the property
$e_0\circ c=c,$ for any $c\in A.$ Take $a:=e_0.$ We have
$$\psi(e_0,b\circ c)=\psi(b,e_0\circ c)=\psi(b,c), \;\; \forall b,c\in A.$$

\noindent Therefore, for $\omega\in C^1_{nov}(A,{\cal K})=C^1(A,{\cal K}),$
such that $\omega(a)=-\psi(e_0,a),$ we have
$$\psi(b,c)=-\omega(b\circ c)=d_{nov}\omega(b,c).$$

\noindent Recall that, $d_{nov}\phi=d_{right}\phi,$ for any
$\phi\in C^1(A,{\cal K}).$ So, any 2-cocycle of $W_1^{nov}$
(in sense of Novikov) with coefficients in the
trivial module, is a coboundary. \quad $\bullet$

\subsection{Cohomologies of $W_n^{rsym}$ in an  antisymmetric module}
Recall that a multiplication in $W_n^{rsym}$ is given by
$a\der_i\circ b\der_j=b\der_j(a)\der_i,$ where
$a,b\in U={\cal K}[[x^{\pm 1},\ldots,x^{\pm 1}]].$ Endow $U$ by a structure of
antisymmetric $W_n^{rsym}$-module: $u\circ a\der_i=a\der_i(u).$

Let $\Omega_n=\{u\,dx_1\wedge\cdots\wedge dx_n: u\in U \}$ be
an antisymmetric $W_n^{rsym}$-module of $n$-dimensional differential forms:
$$(u\,dx_1\wedge\cdots dx_n)\circ a\der_i=
\der_i(au)\,dx_1\wedge\cdots\wedge dx_n.$$

Let $M$ be $W_n^{rsym}$-module.
Construct a cup product of $W_n^{rsym}$-modules
\begin{equation}U\times \Omega^n\otimes M\rightarrow {\bar M}, \;\;
(u\cup v\,dx_1\wedge\cdots\wedge dx_n \,\otimes m)=\pi(uv)m,\label{cupw}
\end{equation}

\noindent where $\pi(u)$ for $u\in U$ denotes a coefficient of $u$ at
$x_1^{-1}\ldots x_n^{-1}.$ Recall that ${\bar M}$ is an antisymmetric
$A$-module corresponding to $M,$ such that ${\bar r}_a=r_a-l_a, {\bar l}_a=0.$
Notice that
$$\pi(der_i(u))=0,\;\; \forall u\in U, i=1,\ldots, n.$$

\noindent Therefore,
$$(u\circ a\der_i)\cup (v\,dx_1\wedge\cdots dx_n\,\otimes m)+
u\cup [v\,dx_1\wedge\cdots dx_n\,\otimes m, a\der_i]=$$
$$a\der_i(u)\cup (v\,dx_1\wedge\cdots\wedge dx_n\otimes m)
+ u\cup (\der_i(av)\,dx_1\wedge\cdots\wedge dx_n\otimes m)$$
$$+u\cup v\,dx_1\wedge\cdots\wedge dx_n\otimes [m,a\der_i]=$$
$$\pi(a\der_i(u)v+u\der_i(av))\,m+\pi(uv)\,[m,a\der_i]=$$
$$\pi(\der_i(auv))+\pi(uv)\,[m,a\der_i]=$$
$$(u\cup v\,dx_1\wedge\cdots\wedge dx_n)\otimes [m,a\der_i]=$$
$$(u\cup v\,dx_1\wedge\cdots\wedge dx_n\otimes m)\circ a\der_i.$$

\noindent and the definition of the cup product~ (\ref{cupw}) is correct.

\begin{th} Let $M$ be an antisymmetric $W_n^{rsym}$-module.
The cup product of $W_n^{rsym}$-modules (\ref{cupw}) induces an isomorphism
$$H^{k+1}_{rsym}(W_n, M)\cong Z^1_{rsym}(W_n,U)\otimes H^k_{lie}(W_n,
\Omega^n\otimes M), \;\; k>0.$$

\noindent An isomorphism takes place
$$Z^1_{rsym}(W_n,U)\cong B^1_{rsym}(W_n,U)=\{dx_i: i=1,\ldots,n\}\cong \wedge ^1.$$
\end{th}

{\bf Proof.} We will argue as in the previos subsection.
For a $W_n$-module $M$ its dual module is denoted by $M'.$
Consider $\wedge^1=\{dx_1,\ldots,dx_n\}$ as a trivial
module over Lie algebra $W_n.$ Endow  $U\otimes \wedge^1$ by a structure of
$W_n$-module using a natural $W_n$-module structure on
$U={\cal K}[[x_1,\ldots,x_n]]$ and a trivial $W_n$-module structure on $\wedge^1:$
$$(u\otimes dx_i)b\der_j=b\der_j(u)\otimes dx_i.$$

It is easy to see that,
$$C^1(W_n,{\cal K})\cong \wedge^1\otimes U',$$

\noindent since for any $f\in C^1(W_n,{\cal K}),$
$$[f,a\der_i](b\der_j)=(f\circ a\der_i)(b\der_j)=-f(a\der_i(b)\der_j).$$

The bilinear map (\ref{cupw}) is nondegenerate and gives a pairing
of modules over Lie algebra $W_n.$ So,
$$U'\cong \Omega_n.$$

Therefore, an isomorphism of $W_n$-modules takes place:
$$C^1(W_n,{\cal K})\cong \wedge^1 \otimes\Omega_n,
\;\; C^1(W_n,M)\cong \wedge^1\otimes \Omega_n\otimes M,$$

\noindent and by theorem~\ref{chev}
$$H^{k+1}_{rsym}(W_n,{M})\cong
\wedge^1\otimes H^k_{lie}(W_n,\Omega^n\otimes M), \;\; k>0.\bullet $$

\begin{crl}
$H^{k+1}_{rsym}(W_n,{\cal K})\cong \wedge^1\otimes
H^k_{lie}(W_n,\Omega^n),\;\; k>0.$ \end{crl}

Recall that $H^*_{lie}(W_n,\Omega^*)$ is calculated by Gelfand
and Fuchs \cite{Fuchs}.

\subsection{Cohomologies of $gl_n^{rsym}$ in an antisymmetric module}
\begin{th} Let $A=gl_n^{rsym}, \; char\,k=0.$ Let  $M$ be a 
finite-dimensional $A$-module, such that a $A^{lie}$-module $M$
is a tensor module. Then the cup product
$M\times {\cal K}\rightarrow M, \;{m\cup\lambda}=\lambda m,$ induces
an isomorphism
$$H^{k+1}_{rsym}(A,M)\cong Z^1_{rsym}(A,M)\otimes H^k_{lie}(A,{\cal K}), \;\; k>0.$$
\end{th}

{\bf Proof.} By theorem~\ref{chev},
$$H^{k+1}_{rsym}(A,M)\cong H^k_{lie}(A,C^1(A,M)), \;\; k>0.$$

By theorem 2.1.2  of \cite{Fuchs},
$$H^k_{lie}(A,C^1(A,M))\cong H^k_{lie}(A,{\cal K})\otimes C^1(A,M)^{A^{lie}}.$$

It remains to notice that,
$$C^1(A,M)^{A^{lie}}:=\{f\in C^1(A,M): [f,a]=0,\;\forall a\in A\}$$

is exactly $Z^1_{rsym}(A,M).$ It is evident:
$$[f,a](b)=(f\circ a)(b)=d_{rsym}f(b,a),\;\; \forall a,b\in A.\quad \bullet$$

\begin{crl} Let $A=gl_n^{rsym}$ and $M$ be a irreducible antisymmetric
$A$-module. Then $H^k_{rsym}(A,M)\ne 0, k\ge 0,$ if and only if $M=A,$  and
$$H^{k+1}_{rsym}(gl_n, (gl_n)_{anti})\cong H^k_{lie}(gl_n,{\cal K}), k\ge 0.$$

\noindent In particular, $H^k_{rsym}(gl_n, {\cal K})=0, \; \; k>0.$
\end{crl}

{\bf Proof.} Since $M$ is antisymmetric,
$$Z^1_{rsym}(A,M)=\{f:A\rightarrow M : f(b\circ a)=f(b)\circ a\}.$$

\noindent Thus, any $f\in Z^1_{rsym}(A,M)$ gives us a homomorphism of right modules
$f:A\rightarrow M.$ Since right modules $M$ and $A$ are irreducible,
by the Lemma of Shur,
$Z^1_{rsym}(A,M)\cong {\cal K},$ if $M\cong A,$ and $Z^1_{rsym}(A,M)=0,$ if
$M\not\cong A.\;\; \bullet$

\subsection{Cohomologies of $gl_n^{rsym}$ in a regular module}
\begin{th} Let $A=gl_n^{rsym}$ over a field ${\cal K}$ of characteristic 0
and $M=A$ be its regular module.
Then the cup product $M\times {\cal K}\rightarrow M, m\cup\lambda=\lambda m,$
induces an isomorphism
$$H^{k+1}_{rsym}(A,M)\cong Z^1_{rsym}(A,A)\otimes H^{k}_{lie}(gl_n,{\cal K}),\;
Z^1_{rsym}(A,A)\cong sl_n, \;\; k\ge 0.$$

\noindent In particular, any cocycle class in $H^{k+1}_{rsym}(gl_n,gl_n)$ has 
a representative that can be presented as $ad\,X\cup \psi,$ where
$\psi\in Z^k_{lie}(gl_n,{\cal K}).$
\end{th}

{\bf Proof.} Any right-symmetric 1-cocycle of an associative algebra $A$
is also an associative 1-cocycle and converse, any associative 1-cocycle
is a right-symmetric 1-cocycle.  So, $Z^1_{rsym}(gl_n,gl_n)\cong Z^1_{ass}(gl_n,gl_n).$
Any derivation of the associative algebra $gl_n$ is a
derivation of the Lie algebra of $gl_n.$ Any derivation of $gl_n^{lie},$
except $a\mapsto tr\,a,$ is inner. So, the following sequence is exact
$$0\rightarrow Z^1_{rsym}(gl_n,gl_n)\rightarrow Z^1_{lie}(gl_n,gl_n)
\rightarrow {\cal K}\rightarrow 0.$$

\noindent In particular,
$$Z^1_{rsym}(gl_n,gl_n)=\{ad\,X: X\in sl_n\}\cong
Z^1_{lie}(sl_n,sl_n)\cong sl_n.$$

It remains to use theorem~\ref{chev}.$\bullet$

\begin{crl} An algebra $gl_n$ as a right-symmetric algebra has
nontrivial deformations. Any right-symmetric local deformation (2-cocycle of
the regular module) is equivalent to a 2-cocycle $\eta_X$ of the form
$$\eta_X(a,b)= (tr\,b)[X,a], $$

\noindent for some $X\in sl_n.$ Any local right-symmetric
deformation can be prolongated. Any formal right-symmetric deformation of
$gl_n$ is equivalent to a deformation of the form
$$\mu_t(a,b)= a\circ b+t\,(tr\,b)[X,a], \quad X\in sl_n,$$

\noindent where $(a,b)\mapsto a\circ b$ is a usual associative
multiplication of matrices.
\end{crl}

{\bf Proof.} These statements can be obtained from the following cohomological
facts
$$H^2_{rsym}(gl_n,gl_n)\cong Z^1_{right}(gl_n,gl_n)\otimes
H^1_{lie}(gl_n,{\cal K}) \cong \{ad\,X\cup tr : X\in sl_n\},$$
$$H^3_{rsym}(gl_n,gl_n)=0,$$

\noindent and corollary~\ref{prod}. \quad $\bullet$

{\bf Remark.} For $\omega\in C^1(gl_n,gl_n),\;\; \omega(a)=(tr\,a)X,$
we have
$$(\eta+d_{rsym}\omega)(a,b)=(tr\,b)[X,a]+(tr\,b)a\circ X
-(tr\,a\circ b)X+(tr\,a)X\circ b=\bar\eta_X(a,b),$$

\noindent where
$$\bar\eta_X(a,b)= (tr\,b)X\circ a-(tr\,a\circ b)X+(tr\,a)X\circ b.$$

\noindent Therefore, $[\eta_X]\sim [\bar\eta_X].$ Notice that $\bar\eta_X$
is a symmetric cocycle:
$$\bar\eta(a,b)=\bar\eta(b,a).$$

\noindent The prolongation formula for $\bar\eta_X$ is a little bit complicated. It can be
given by
$$\bar\mu_t(a,b)=\Phi_t^{-1}(\mu_t(\Phi_t(a),\Phi_t(b))),$$

\noindent where
$$\Phi_t=id+t\,\omega.$$

\noindent In particular,
$$\Phi^{-1}_t(a)=a-t\,(tr\,a)X+t^2\,(tr\,a)^2X^2-\cdots,$$

\noindent and some of the beginning terms of $\bar\mu_t$ look like
$$\bar\mu_t(a,b)=a\circ b+t\,X\circ((tr\,a)b-(tr\,a\circ b)+(tr\,b)a)$$
$$+t^2\,(tr\,a\,tr\,b-(tr\,a\circ b)^2X^2-(tr\,a\,tr\,(X\circ b))X
-(tr\,(a\circ X)\,tr\,b)X)+\cdots $$

So, right-symmetric prolongation can be constructed in a such way, that
corresponding Lie multiplication will not be changed:
$$\bar\mu_t(a,b)-\bar\mu_t(b,a)=[a,b].$$

{\bf Acknowledgements.} {\sl The author is deeply grateful to INTAS foundation 
(prof. U.Rehmann), Bielelfeld University  (prof. K.Ringel) and the Abdus Salam International Centre 
for Theoretical Physics (prof.M.Virasoro, prof. M.S. Narasimhan) 
for support.}

\end{document}